\documentclass{amsart}

\input{prepictex}
\input{pictex}
\input{postpictex}

\begin{document}

\def\sbullet{\raise2.3pt\hbox{$_\bullet$}}
\def\bddots{\kern-10pt\lower7pt\hbox{.}\kern-7pt\lower11pt\hbox{.}
\kern-7pt\lower15pt\hbox{.}}

\def\CC{\mathbb C}
\def\NN{\mathbb N}
\def\RR{\mathbb R}
\def\TT{\mathbb T}
\def\ZZ{\mathbb Z}

\def\Bb{{\mathcal B}}
\def\Mm{{\mathcal M}}
\def\Ff{{\mathcal F}}
\def\Hh{{\mathcal H}}
\def\Kk{{\mathcal K}}
\def\Oo{{\mathcal O}}
\def\Ss{{\mathcal S}}
\def\Tt{{\mathcal T}}
\def\Ad{\operatorname{Ad}}
\def\Aut{\operatorname{Aut}}
\def\clsp{\overline{\lsp}}
\def\End{\operatorname{End}}
\def\id{\operatorname{id}}
\def\Isom{\operatorname{Isom}}
\def\Ker{\operatorname{Ker}}
\def\lsp{\operatorname{span}}
\def\Hom{\operatorname{Hom}}
\def\Obj{\operatorname{Obj}}
\def\MCE{\operatorname{MCE}}
\def\fin{\operatorname{fin}}
\def\Ext{\operatorname{Ext}}
\def\L{\Lambda}
\def\uL(#1){\Lambda^{\le #1}}
\def\Lmin(#1,#2){{%
\Lambda^{\min}(#1,#2) }}
\def\MU(#1,#2;#3,#4){{%
\Theta(#1)^{\Pi #2}_{#3,#4} }}%

\theoremstyle{plain}
\newtheorem{theorem}{Theorem}[section]
\newtheorem*{theorem*}{Theorem}
\newtheorem*{prop*}{Proposition}
\newtheorem{cor}[theorem]{Corollary}
\newtheorem{lemma}[theorem]{Lemma}
\newtheorem{prop}[theorem]{Proposition}
\theoremstyle{remark}
\newtheorem{rmk}[theorem]{Remark}
\newtheorem{rmks}[theorem]{Remarks}
\newtheorem*{aside}{Aside}
\newtheorem*{note}{Note}
\newtheorem{comment}[theorem]{Comment}
\newtheorem{example}[theorem]{Example}
\newtheorem*{example*}{Example}
\newtheorem{examples}[theorem]{Examples}
\theoremstyle{definition}
\newtheorem{dfn}[theorem]{Definition}
\newtheorem{dfns}[theorem]{Definitions}
\newtheorem{notation}[theorem]{Notation}
\numberwithin{equation}{section}

\title[${C^*}$-algebras of higher-rank graphs] {\boldmath{The 
${C^*}$-algebras of finitely aligned higher-rank graphs}}

\author{Iain Raeburn}
\author{Aidan Sims}
\author{Trent Yeend}

\address{Department of Mathematics  \\ University of Newcastle\\  NSW 
2308\\ AUSTRALIA}
\email{iain,aidan,trent@maths.newcastle.edu.au}
\date{May 26, 2003}
\subjclass{Primary 46L05}
\thanks{This research was supported by the Australian Research Council.}
\keywords{graph algebra; Cuntz-Krieger algebra; uniqueness}

\begin{abstract} We generalise the theory of Cuntz-Krieger families
and graph algebras to the class of \emph{finitely aligned} $k$-graphs. 
This class contains in particular all row-finite $k$-graphs. The Cuntz-Krieger
relations for non-row-finite $k$-graphs look significantly different from
the usual ones, and this substantially complicates
the analysis of the graph algebra. We prove 
a gauge-invariant uniqueness theorem and a  
Cuntz-Krieger uniqueness  theorem for the $C^*$-algebras of 
finitely  aligned $k$-graphs.
\end{abstract}

\maketitle

\section{Introduction}

It has been known for many years that the Cuntz-Krieger algebras 
of (0,1)-matrices
\cite{CK} can be viewed as the $C^*$-algebras of directed graphs 
\cite{EW}. More
recently, the construction has been extended to cover infinite directed graphs
\cite{KPR, FLR} and higher-rank analogues, known as $k$-graphs \cite{KP}. The
resulting classes of \emph{graph algebras} contain many interesting 
examples, and have
in particular provided a rich supply of models for the classification 
theory of simple
purely infinite nuclear
$C^*$-algebras \cite{Sz}.

Graph algebras have now been associated to all infinite graphs, and an elegant
structure theory relates the behaviour of loops in a graph to the
properties of its graph algebra. For $k$-graphs, the current state of affairs is less
satisfactory. The object of this paper is to associate graph algebras 
to a wide class
of infinite $k$-graphs, and to prove versions of the gauge-invariant uniqueness
theorem and the Cuntz-Krieger uniqueness theorem for these graph algebras.

Before describing our approach, we recall how the theory of graph 
algebras developed.
A directed graph
$E$ consists of a countable vertex set
$E^0$, a countable edge set
$E^1$, and range and source maps $r,s : E^1 \to E^0$. When each vertex receives
at most finitely many edges ($E$ is \emph{row-finite})
the graph algebra $C^*(E)$ is the universal $C^*$-algebra generated by mutually
orthogonal projections $\{p_v : v
\in E^0\}$ and partial isometries $\{s_e : e \in
E^1\}$ satisfying $s^*_e s_e = p_{s(e)}$ for all $e \in E^1$ and
\begin{equation}\label{eqn:CKrel} p_v = \sum_{r(e) = v} s_e s^*_e\ \mbox{ when
$r^{-1}(v)$ is non-empty.}
\end{equation}
When $r^{-1}(v)$ is infinite, the sum on the 
right-hand side of
\eqref{eqn:CKrel} cannot converge in a $C^*$-algebra, and hence the relation
must be adjusted.  The appropriate adjustment 
was suggested by the
analysis of the Toeplitz algebras of Hilbert bimodules in
\cite{FR}: impose relation \eqref{eqn:CKrel} only where
$r^{-1}(v)$ is finite, and add the requirement that the $s_e$ have 
orthogonal range
projections dominated by $p_{r(e)}$ (which in the row-finite case follows from
\eqref{eqn:CKrel}). The resulting family of graph algebras was studied in
\cite{FLR}. That these are the appropriate relations was confirmed when other
authors with different points of view arrived at the same conclusion 
\cite{Pat, Spiel}.

The first work on higher-rank graphs concerned row-finite $k$-graphs without
sources \cite{KP}. For directed graphs (that is, when $k=1$), there 
is a constructive procedure
for extending results to graphs with sources \cite[Lemma~1.2]{BPRS}. 
However when $k>1$, there
are many different kinds of sources, and there is as yet no analogous 
procedure for dealing
with them. In \cite{RSY1}, we considered a class of row-finite
$k$-graphs which may have sources provided a local convexity 
condition is satisfied.  In
\cite{RS1}, Raeburn and Sims studied infinite $k$-graphs by viewing 
them as product
systems of graphs, as in \cite{FS}, and applying the 
techniques of \cite{F99}
to the Toeplitz algebras of the associated product system of Hilbert 
bimodules. The
analysis in \cite{RS1} led to two conclusions. First, it identified an
extra Cuntz-Krieger relation which is automatic for row-finite 
$k$-graphs, but is not in
general. This extra relation is needed to ensure that the algebras generated by
Cuntz-Krieger families are spanned by partial isometries of the usual form.
Unfortunately, the new relation can involve infinite sums of 
projections (see
\cite[Remark~7.2]{RS1}); the second conclusion of
\cite{RS1} was that we should restrict attention to the \emph{finitely aligned}
$k$-graphs for which the new relation is
$C^*$-algebraic rather than spatial.

In this paper we introduce Cuntz-Krieger relations which are appropriate for
arbitrary finitely aligned $k$-graphs. We do not assume that our
$k$-graphs are locally convex or row-finite, and we do allow them to 
have  sources. When
$k=1$ or the
$k$-graph is row-finite and locally convex, our new Cuntz-Krieger relations are
equivalent to the usual ones. We show that for every finitely
aligned $k$-graph $\L$, there is a family of nonzero partial isometries
which satisfies the new relations, and we define $C^*(\L)$ to be the
universal $C^*$-algebra generated by such a family.  We then prove 
versions of the
gauge-invariant uniqueness theorem and the Cuntz-Krieger uniqueness theorem for
$C^*(\L)$. Our analysis is elementary in the sense that we do not use 
groupoids, partial
actions or Hilbert bimodules, though we cheerfully acknowledge that 
we have gained
insight from the models these theories provide.

The results in this paper extend the existing theory of graph 
algebras in several
directions. Since
$1$-graphs are always finitely aligned, and our new relations are 
then equivalent to the
usual ones (Proposition~\ref{prop:new rels
generalise}), our approach provides the first
elementary analysis of the
$C^*$-algebra of an arbitrary directed graph. Our results are also 
new for finitely
aligned
$k$-graphs without sources; those interested primarily in this 
situation may mentally
replace all the symbols $\L^{\leq n}$ by $\L^n$, and thereby avoid 
several technical
complications. Even for row-finite
$k$-graphs we make significant improvements on the existing theory: for
non-locally-convex row-finite $k$-graphs, our Cuntz-Krieger families 
may have every
vertex projection nonzero, unlike those in \cite{RSY1} (see 
Example~\ref{ex:non-loc-conv}).

\smallskip

In Section~\ref{sec:kgraphs and repns} we describe our new 
Cuntz-Krieger relations for
a finitely aligned
$k$-graph $\L$, define $C^*(\L)$ to be
the universal $C^*$-algebra generated by a Cuntz-Krieger family, and 
investigate some of
its basic properties. We discuss a notion of boundary paths which we 
use to construct a
Cuntz-Krieger family in which every vertex projection is nonzero.

The \emph{core} in $C^*(\L)$ is the fixed-point algebra 
$C^*(\L)^\gamma$ for
the gauge action $\gamma$ of $\TT^k$. In Section~\ref{sec:thecore} we 
show that the
core is AF, and deduce that a homomorphism
$\pi$ of
$C^*(\L)$ which is nonzero at each vertex projection is injective on the core.

Our proof that $C^*(\L)^\gamma$ is AF is quite different from the 
argument which we gave
for row-finite $k$-graphs in \cite{RSY1} in that we do not describe
$C^*(\L)^\gamma$ as a direct limit over $\NN^k$. Instead, we describe 
$C^*(\L)^\gamma$ as
the increasing union of finite-dimensional algebras indexed by finite sets of paths,
and produce families of matrix units which span these algebras. In addition to 
showing that $C^*(\L)$ is AF, this formulation is a key ingredient in our proof of the Cuntz-Krieger uniqueness theorem. The uniqueness
theorems themselves are proved in
Section~\ref{sec:theorems}.

We conclude with three appendices in which we discuss various aspects 
of our new
Cuntz-Krieger relations.
In Appendix~\ref{app:why new rel} we explain our motivation for 
introducing these new
  and apparently substantially different relations; we describe 
examples illustrating
the other possibilities we considered, and their failings. In 
Appendix~\ref{app:new rels
generalise}, we show that for ordinary directed graphs (that is, for $k=1$) and
for locally convex row-finite $k$-graphs, our new Cuntz-Krieger 
relations are equivalent
to the usual ones. Appendix~\ref{app:gens do} gives an equivalent 
formulation of our
Cuntz-Krieger relations using only the edges in the $1$-skeleton of 
the $k$-graph.

\section{$k$-graphs and Cuntz-Krieger families} \label{sec:kgraphs and 
repns} We regard
$\NN^k$ as a semigroup with identity $0$. For $1 \le i \le k$, we 
write $e_i$ for the
$i^{\rm{th}}$ generator of $\NN^k$, and for $n \in \NN^k$ we write 
$n_i$ for the
$i^{\rm{th}}$ coordinate of $n$. We use $\le$ for the partial order 
on $\NN^k$ given by
$m \le n$ if $m_i \le n_i$ for all $i$. The expression $m < n$ means 
$m \le n$ and $m
\not= n$, and does not necessarily indicate that $m_i < n_i$ for all 
$i$. For $m,n \in
\NN^k$, we write $m \vee n$ for their coordinate-wise maximum and $m 
\wedge n$ for their
coordinate-wise minimum.

A \emph{$k$-graph} is a pair $(\L, d)$ consisting of a countable 
small category $\L$ and
a \emph{degree functor} $d: \L \to \NN^k$ which satisfy the \emph{factorisation
property:} for every $\lambda \in \L$ and $m,n \in \NN^k$ with 
$d(\lambda) = m + n$
there exist unique $\mu,\sigma \in \L$ such that $d(\mu) = m$, 
$d(\sigma) = n$ and
$\lambda = \mu\sigma$. 

Since we are regarding $\L$ as a type of 
graph, we refer to the
morphisms of $\L$ as paths and to the objects of $\L$ as 
vertices, and write $s$
and $r$ for the domain and codomain maps. For a thorough introduction 
to the structure
of $k$-graphs, see \cite[Section~2]{RSY1}.

\begin{notation} We use lower-case Greek letters to denote paths in
$k$-graphs. However, we reserve $\delta$ for the Kronecker delta, and 
$\gamma$ for the gauge action (see Section~\ref{sec:thecore}).
\end{notation}

Given $k$-graphs $(\L,d_\L)$ and $(\Gamma,d_\Gamma)$, a \emph{graph 
morphism} from $\L$
to $\Gamma$ is a functor $x : \L \to \Gamma$ such that $d_\Gamma(x(\lambda)) =
d_\L(\lambda)$ for all $\lambda \in \L$.    For $n \in \NN^k$, $\L^n$ is the 
collection of all paths of degree $n$; that is
\[
\L^n := \{\lambda \in \L : d(\lambda) = n\}.
\]

The factorisation property ensures that associated to each vertex $v 
\in \Obj(\L)$ there
is a unique element of $\L^0$ whose range (and hence source) is $v$; 
we call this
morphism $v$ as well, identifying $\Obj(\L)$ with $\L^0$.  For $E 
\subset \L$ and
$\lambda \in \L$, we define
\begin{gather*}
\lambda E :=\{\lambda\mu : \mu \in E, r(\mu) = s(\lambda)\},\text{ 
and } \\ E \lambda
:=\{\mu\lambda : \mu \in E, s(\mu) = r(\lambda)\}.
\end{gather*} Hence, for $v \in \L^0$ and $E \subset \L$, $v E = 
\{\mu \in E : r(\mu) =
v\}$ and $E v = \{\mu \in E : s(\mu) = v\}$.

For $n \in \NN^k$, we define
\[
\uL(n) := \{\lambda \in \L : d(\lambda) \le n,\text{ and } 
d(\lambda)_i < n_i \implies
s(\lambda)\L^{e_i} = \emptyset\}.
\]

For $\lambda \in \L$ and $m \le n \le d(\lambda)$, the factorisation 
property gives
unique paths $\lambda' \in \L^m$, $\lambda'' \in \L^{n-m}$ and $\lambda''' \in
\L^{d(\lambda) - n}$ such that $\lambda = 
\lambda'\lambda''\lambda'''$. We denote
$\lambda''$ by $\lambda(m,n)$, so $\lambda' = \lambda(0,m)$ and $\lambda''' =
\lambda(n,d(\lambda))$. More generally, for all
$m \le n \in \NN^k$, $\lambda(m,n) := \lambda(m
\wedge d(\lambda), n \wedge d(\lambda))$.

\begin{dfn}\label{dfn:Lmin} For $\lambda,\mu \in \L$, we write
\[
\Lmin(\lambda,\mu) := \{(\alpha,\beta) : \lambda\alpha = \mu\beta, 
d(\lambda\alpha) =
d(\lambda) \vee d(\mu)\}
\] for the collection of pairs which give \emph{minimal common 
extensions} of $\lambda$
and $\mu$. We say that $\L$ is \emph{finitely aligned} if 
$\Lmin(\lambda,\mu)$ is finite
(possibly empty) for all $\lambda,\mu \in \L$.
\end{dfn}

\begin{rmk} For $\lambda,\mu \in \L$, the map $(\alpha,\beta) \mapsto 
\lambda\alpha$ is
a bijection between $\Lmin(\lambda,\mu)$ and the set 
$\MCE(\lambda,\mu)$ defined in
\cite[Definition~5.3]{RS1}. Hence our definition of a finitely aligned $k$-graph 
agrees with that of
\cite{RS1}.
\end{rmk}

\begin{dfn}\label{dfn:exhaustive} Let $(\L,d)$ be a $k$-graph, let $v 
\in \L^0$ and $E
\subset v\L$. We say that $E$ is \emph{exhaustive} if for every $\mu 
\in v\L$ there
exists $\lambda \in E$ such that $\Lmin(\lambda,\mu)\not= \emptyset$.
\end{dfn}

\begin{dfn}\label{dfn:CKfamily} Let $(\L,d)$ be a finitely aligned $k$-graph. A
\emph{Cuntz-Krieger $\L$-family} is a 
collection $\{t_\lambda :
\lambda \in \L\}$ of partial isometries in a $C^*$-algebra satisfying
\begin{itemize}
\item[(i)] $\{t_v : v \in \L^0\}$ is a collection of mutually 
orthogonal projections;
\item[(ii)] $t_\lambda t_\mu = t_{\lambda\mu}$ whenever $s(\lambda)=r(\mu)$;
\item[(iii)] $t^*_\lambda t_\mu = \sum_{(\alpha,\beta) \in 
\Lmin(\lambda,\mu)} t_\alpha
t^*_\beta$ for all $\lambda,\mu \in \Lambda$; and
\item[(iv)] $\prod_{\lambda \in E} (t_v - t_\lambda t^*_\lambda) = 0$ 
for all $v \in
\L^0$ and finite exhaustive $E \subset v\L$.
\end{itemize}
\end{dfn}

\begin{rmk}\label{rmk:why new rels} A number of aspects of these 
Cuntz-Krieger relations
are worth commenting on:
\begin{itemize}
\item As seen in \cite{RS1}, the restriction to finitely aligned 
$k$-graphs is necessary
for the sum in relation (iii) to make sense.
\item Relation (iii) implies that 
$t^*_\lambda t_\lambda = t_{s(\lambda)}$, and that
$t^*_\lambda t_\mu = 0$ if $\Lmin(\lambda,\mu) = \emptyset$.
\item Relations (iii) and (iv) have been significantly changed from 
their usual form
(see \cite[Section~1]{BPRS} and \cite[Definition~3.3]{RSY1}), and we 
feel they require
explanation. The short explanation is that they are the right 
relations for generating
tractable Cuntz-Krieger algebras for which a homomorphism is 
injective on the core if
and only if it is nonzero at each vertex projection 
(Theorem~\ref{thm:faithful on core}). A much more 
detailed explanation is
contained in Appendix~\ref{app:why new rel}.
\item In Appendix~\ref{app:new rels generalise} we prove that for 
$1$-graphs and for
locally convex row-finite $k$-graphs, our relations are equivalent to 
those set forth in
\cite{FLR} and \cite{RSY1} respectively.
\item Previous treatments of $k$-graph $C^*$-algebras have shown that 
the Cuntz-Krieger
relations can be formulated in terms of the $1$-skeleton of $\L$; that is
in terms of vertices and paths of degree $e_i$. 
We show in Appendix~\ref{app:gens do} that the 
same is true for our relations.
\end{itemize}
\end{rmk}

Given a finitely aligned $k$-graph $(\L,d)$, there exists a 
$C^*$-algebra $C^*(\L)$
generated by a Cuntz-Krieger $\L$-family $\{s_\lambda : \lambda \in 
\Lambda\}$ which is
universal in the following sense: given a Cuntz-Krieger $\L$-family 
$\{t_\lambda :
\lambda \in \Lambda\}$, there exists a unique homomorphism $\pi_t$ 
of $C^*(\L)$ such that
$\pi_t(s_\lambda) = t_\lambda$ for all $\lambda \in \L$.

The following lemma sets forth some useful consequences of
Definition~\ref{dfn:CKfamily}(i)--(iii).

\begin{lemma}\label{lem:uptoLambda orth range projections} Let 
$(\L,d)$ be a finitely
aligned $k$-graph and let $\{t_\lambda:\lambda\in\L\}$ be a family of 
partial isometries
satisfying Definition~\ref{dfn:CKfamily}(i)--(iii). Then
\begin{itemize}
\item[(i)]
$t_\lambda t^*_\lambda t_\mu t^*_\mu=\sum_{(\alpha,\beta)\in\Lmin(\lambda,\mu)}
t_{\lambda\alpha}t^*_{\lambda\alpha}$ for all $\lambda,\mu\in\L$. In 
particular,
$\{t_\lambda t^*_\lambda :\lambda\in\L\}$ is a family of commuting projections.
\item[(ii)] For $\lambda,\mu\in\uL(n)$, we have $t^*_\lambda
t_\mu=\delta_{\lambda,\mu}t_{s(\lambda)}$.
\item[(iii)] If $E\subset v\uL(n)$ is finite, then $t_v\ge 
\sum_{\lambda\in E}t_\lambda
t^*_\lambda$.
\item[(iv)]
$C^*(\{t_\lambda : \lambda \in \L\})
=\clsp\{t_\lambda t^*_\mu :\lambda,\mu\in\L\}
=\clsp\{t_\lambda t^*_\mu :\lambda,\mu\in\L, s(\lambda) = s(\mu)\}$.
\end{itemize}
\end{lemma}

\begin{proof} Part (i) is obtained by multiplying both sides of the equation in
Definition~\ref{dfn:CKfamily}(iii) on the left by $t_\lambda$ and on 
the right by
$t^*_\mu$.

For (ii), suppose that $t^*_\lambda t_\mu \not= 0$. Then
Definition~\ref{dfn:CKfamily}(iii) ensures that there exists 
$(\alpha,\beta) \in
\Lmin(\lambda,\mu)$, so $\lambda\alpha = \mu\beta$ and 
$d(\lambda\alpha) \le n$. Since
$\lambda,\mu \in \uL(n)$, it follows that $\alpha = \beta = s(\lambda)$, so
$\lambda=\mu$.

For (iii), note that if $\lambda,\mu \in E$ and $\lambda \not= \mu$, 
then $t_\lambda
t^*_\lambda t_\mu t^*_\mu = 0$ by (ii), and $t_v t_\lambda 
t^*_\lambda = t_\lambda
t^*_\lambda$ for all $\lambda \in E$ by Definition~\ref{dfn:CKfamily}(ii).

For part (iv), note that $\clsp\{t_\lambda t^*_\mu :\lambda,\mu\in\L\}$ is 
clearly closed under
adjoints and contains $\{t_\lambda : \lambda \in \L\}$. Furthermore,
$\clsp\{t_\lambda : \lambda \in \L\}$  is
closed under multiplication by
Definition~\ref{dfn:CKfamily}(iii). To see that
$\clsp\{t_\lambda t^*_\mu :\lambda,\mu\in\L\}
=\clsp\{t_\lambda t^*_\mu :\lambda,\mu\in\L, s(\lambda) = s(\mu)\}$, note
that if $s(\lambda) \not= s(\mu)$ then $t_\lambda t^*_\mu = t_\lambda 
t_{s(\lambda)} t_{s(\mu)}^* t_\mu^* = 0$ by Definition~\ref{dfn:CKfamily}(i).
\end{proof}

We define our prototypical Cuntz-Krieger $\L$-family using a 
boundary-path space
associated to $\L$. For $m \in (\NN \cup\{\infty\})^k$, 
recall from \cite[Examples~2.2(ii)]{RSY1} the definition of the
$k$-graph $\Omega_{k,m}$:
\begin{gather*}
\Obj(\Omega_{k,m}) = \{p \in \NN^k : p \le m\}, \\
\Hom(\Omega_{k,m}) = \{(p,q) \in \Obj(\Omega_{k,m}) \times 
\Obj(\Omega_{k,m}) : p \le q\}, \\ 
r(p,q) = p, \quad s(p,q) = q, \quad d(p,q) = q-p.
\end{gather*}

If $x : \Omega_{k,m} \to \L$ is a graph morphism and $\lambda \in \L$ 
with $s(\lambda) = x(0)$, then there is a unique graph morphism
$\lambda x : \Omega_{k, m+d(\lambda)} \to \L$ such that 
$(\lambda x)(0,d(\lambda)) = \lambda$, and $(\lambda  x)(d(\lambda),n) 
= x(0, n-d(\lambda))$ for all $n \ge d(\lambda)$. If $x : \Omega_{k,m} \to 
\L$ is a graph morphism and $n \in \NN^k$ with $n \le m$, then there is a 
unique graph morphism $x(n,m) : \Omega_{k,m-n} \to \L$ such that 
$\big(x(n,m)\big)(0,l) = x(n, n+l)$ for all $l \in \NN^k$. Notice that these
two constructions are inverse in the sense that
$(\lambda x)\big(d(\lambda), d(\lambda x)\big)$ and $x(0,n) x(n, m)$ are 
both equal to $x$.

\begin{dfn}\label{dfn:bdry path} Let $(\L,d)$ be a $k$-graph, let $m \in
(\NN\cup\{\infty\})^k$, and let $x : \Omega_{k,m} \to \L$ be a graph 
morphism. We call
$x$ a \emph{boundary path} if there exists $n_x \in \NN^k$ such that 
$n_x \le m$ and
\begin{equation}\label{eqn:bdry path} n \in \NN^k, n_x \le n \le m 
\text{ and } n_i =
m_i \text{ imply that } x(n)\L^{e_i} = \emptyset.
\end{equation} We extend the range and degree maps to boundary paths 
$x : \Omega_{k,m}
\to \L$ by setting $r(x) := x(0)$ and $d(x) := m$. We write 
$\uL(\infty)$ for the
collection of all boundary paths of $\Lambda$, and $v\uL(\infty)$ for $\{x \in
\uL(\infty) : r(x) = v\}$.
\end{dfn}

\begin{rmk}
If $\L$ has no sources, then the boundary path space $\uL(\infty)$ is the usual 
infinite path space $\L^\infty$ of \cite[Definitions~2.1]{KP} consisting of all graph 
morphisms $x : \Omega_{k,(\infty,\dots,\infty)}  \to \L$.
\end{rmk}

\begin{lemma}
Let $(\L,d)$ be a $k$-graph, and let $x \in \uL(\infty)$. 
\begin{itemize}
\item[(i)] If $\lambda \in \L$ with $s(\lambda) = r(x)$, then $\lambda x \in \uL(\infty)$.
\item[(ii)]  If $n \in \NN^k$ with $n \le d(x)$, then $x(n,d(x)) \in \uL(\infty)$.
\end{itemize}
\end{lemma}
\begin{proof}
We need only show that there exist $n_{\lambda x}$ and $n_{x(n,d(x))}$
satisfying \eqref{eqn:bdry path}. This works with $n_{\lambda x} := n_x 
+ d(\lambda)$ and $n_{x(n,d(x))} := (n_x - n) \vee 0$.
\end{proof}

\begin{lemma}\label{lem:uLinfty(v) nonempty} 
Let $(\L,d)$ be a $k$-graph. Then $v\uL(\infty)$ is nonempty for all $v \in \L^0$.
\end{lemma}

\begin{proof} 
For $i \in \NN$ write $[i]$ for the element of 
$\{1,\dots,k\}$ which is
congruent to $i$ (mod $k$). Fix $v \in \L^0$. Construct a sequence of paths with range 
$v$ as follows:
$\lambda_0 := v$, and given $\lambda_{i-1}$,
\[
\lambda_{i} := \lambda_{i-1}\nu\text{ for some $\nu \in 
s(\lambda_{i-1})\uL(e_{[i]})$;}
\] so at the $i^{\rm{th}}$ step, we append a segment of degree 
$e_{[i]}$ if possible,
and append nothing otherwise. 

Define $m := \lim_{i \to \infty} d(\lambda_i) \in (\NN \cup\{\infty\})^k$. Then 
there is a unique graph morphism $x : \Omega_{k,m} \to \L$ such that $x(0,
d(\lambda_i)) = \lambda_i$ for all $i \in \NN$. To show that $x$ is a 
boundary path, we need only produce $n_x \in \NN^k$ with $n_x \le m$ 
which  satisfies \eqref{eqn:bdry path}.

For each $j\in\{1,\dots,k\}$ such that $s(\lambda_{i-1})\L^{e_j}=\emptyset$ for some $i$, let
\[
 i(j) := \min\{i \in \NN : [i] = j\text{ and } s(\lambda_{i-1})\L^{e_j} = \emptyset\}.
\] 
Let $I := \max\{i(j) : m_j < \infty\}$, and let $n_x 
:= d(\lambda_I)$.

Suppose that $n \in \NN^k$ with $n_x \le n \le m$, and that $n_j = 
m_j$. Then $m_j <
\infty$ so $i(j)$ is defined and $I \ge i(j)$ by definition. Since $n \ge n_x =
d(\lambda_I)$, it follows that $n \ge d(\lambda_{i(j) - 1})$. But
$s(\lambda_{i(j)-1})\L^{e_j}= \emptyset$, which implies $x(n)\L^{e_j} 
= \emptyset$ by
the factorisation property.
\end{proof}

\begin{prop}\label{prp:bdry repn nonzero} Let $(\L, d)$ be a finitely 
aligned $k$-graph.
For $\lambda \in \L$, define
\[
 S_\lambda e_x := 
\begin{cases} 
e_{\lambda x} &\text{if $s(\lambda)=r(x)$} \\ 
0 &\text{otherwise.}
\end{cases}
\]
Then $\{S_\lambda : \lambda \in \L\}$ 
is a Cuntz-Krieger $\L$-family called the \emph{boundary-path
representation}. Furthermore, every $S_v$ is nonzero.
\end{prop}
\begin{proof} 
It follows from Lemma~\ref{lem:uLinfty(v) nonempty} 
that each $S_v$ is nonzero. 

A simple calculation using inner products in $\ell^2(\uL(\infty))$ shows that
\[
S^*_\lambda e_x = 
\begin{cases}
e_{x(d(\lambda), d(x))} &\text{if $x(0, d(\lambda)) = \lambda$} \\
0 &\text{otherwise.}
\end{cases}
\]

We need to check (i)--(iv) of Definition~\ref{dfn:CKfamily}.

Relation (i) holds since $S_v$ is the projection onto $\clsp\{e_x : x \in 
v\uL(\infty)\}$.

Checking (ii) amounts to showing that the boundary path $\lambda(\mu 
x)$ is equal to the
boundary path $(\lambda\mu) x$. This follows from associativity of 
composition in the category $\L$.

Relation (iii) follows from a simple calculation involving inner products (see
\cite[Example~7.4]{RS1}).

To check that (iv) holds, let $E \subset v\L$ be finite and exhaustive and 
let $x \in
v\uL(\infty)$. It suffices to show that $\prod_{\lambda \in E} (S_v - S_\lambda
S^*_\lambda) e_x = 0$. Let 
\[ \textstyle
N :=\big(\bigvee_{\lambda\in E}d(\lambda)\big) 
\vee n_x; 
\]
in particular, $N \ge n_x$ so \eqref{eqn:bdry path} implies
$x(N) \L^{e_j} = \emptyset$ whenever $m_j < \infty$.  
Since $E$ is exhaustive, there exists $\lambda_x\in E$ such that
$\Lmin(x(0,N),\lambda_x)\neq\emptyset$; let $(\alpha,\beta) \in 
\Lmin(x(0,N),\lambda_x)$. We claim that $\alpha = x(N)$. Suppose for 
contradiction $d(\alpha)_i > 0 $ for some $i$. Then 
$d(x(0,N))_i < d(\lambda_x)_i$. But $N_i \ge d(\lambda_x)_i$ by definition,
and hence we must have $d(x)_i < N_i$, so $m_i < \infty$. 
Hence $x(N)\L^{e_i} = \emptyset$ contradicting 
$d(\alpha)_i > 0$. This establishes the claim, giving $x(0,N) = \lambda_x\beta$,
and hence $x(0, d(\lambda_x)) = \lambda_x$. But then
\[
\bigg(\prod_{\lambda \in E} (S_v - S_\lambda S^*_\lambda)\bigg)e_x =
\bigg(\prod_{\lambda \in E\setminus\{\lambda_x\}} (S_v - S_\lambda 
S^*_\lambda)\bigg)
(S_v - S_{\lambda_x} S^*_{\lambda_x}) e_x = 0
\] 
because  $S_v e_x = e_x = S_{\lambda_x} S^*_{\lambda_x} e_x$.
\end{proof}

\section{Analysis of the core} \label{sec:thecore}

Given a finitely aligned $k$-graph $(\L,d)$, there is a strongly 
continuous \emph{gauge
action} $\gamma:\TT^k\to\Aut(C^*(\L))$ determined by $\gamma_z(s_\lambda) =
z^{d(\lambda)}s_\lambda$ where $z^m = z_1^{m_1}\cdots z_k^{m_k} \in \TT$. The 
fixed-point
algebra $C^*(\L)^\gamma$ is equal to $\clsp\{s_\lambda s^*_\mu : 
d(\lambda) = d(\mu)\}$
and is called the \emph{core} of $C^*(\L)$.

\begin{theorem}\label{thm:faithful on core} Let $(\L,d)$ be a finitely 
aligned $k$-graph.
Then $C^*(\L)^\gamma$ is AF. If $\{t_\lambda : \lambda \in \Lambda\}$ 
is a Cuntz-Krieger
$\L$-family with $t_v \not= 0$ for all $v \in \L^0$, then the 
homomorphism $\pi_t$ of
$C^*(\L)$ such that $\pi_t(s_\lambda) = t_\lambda$ is injective on 
$C^*(\L)^\gamma$.
\end{theorem}

The remainder of this section 
is devoted to proving Theorem~\ref{thm:faithful on core}. We
therefore fix a finitely aligned $k$-graph $(\L,d)$ and a Cuntz-Krieger $\L$-family 
$\{t_\lambda : \lambda \in \L\}$. We also fix a finite set $E \subset \L$. 
We want to identify a finite set $\Pi E$ containing
$E$ such that $\lsp\{s_\lambda s^*_\mu:\lambda,\mu\in\Pi E, 
d(\lambda)=d(\mu)\}$ is
closed under multiplication, and hence is a finite-dimensional subalgebra of
$C^*(\L)^\gamma$. The next Lemma implies that such sets exist.

\begin{lemma}\label{lem:Pi E} There exists a finite set $F \subset \L$ 
which contains $E$ and satisfies
\begin{equation}\label{eqn:Pi E}
\begin{split}
\lambda,\mu&,\sigma\!,\tau \in F, d(\lambda) = d(\mu), d(\sigma) = 
d(\tau), s(\lambda)  =
s(\mu) \\ &\text{ and } s(\sigma) = s(\tau) \text{ imply } 
\{\lambda\alpha, \tau\beta :
(\alpha,\beta) \in \Lmin(\mu,\sigma)\}\subset F.
\end{split}
\end{equation}
Moreover, for any finite $F$ which contains $E$ and satisfies \eqref{eqn:Pi E},
\[
M^t_F := \lsp\{t_\lambda t^*_\mu : \lambda,\mu 
\in F, d(\lambda) = d(\mu)\}
\]
is a finite-dimensional $C^*$-subalgebra of 
$C^*(\{t_\lambda t^*_\mu : d(\lambda) = d(\mu)\})$.
\end{lemma}

Before proving Lemma~\ref{lem:Pi E}, we recall from 
\cite[Definition~8.3]{RS1} that for $F \subset \L$, 
\[
\MCE(F) := \{\lambda \in \L : d(\lambda) = \textstyle{\bigvee_{\alpha
\in F}}\,d(\alpha) \text{ and } \lambda(0, d(\alpha)) = 
\alpha\text{ for all } \alpha
\in F\},
\] 
and that $\vee F := \bigcup_{G \subset F} \MCE(G)$. Lemma~8.4 of 
\cite{RS1} shows that $\vee F$ contains $F$, 
is finite whenever $F$ is, and is
closed under taking minimal common extensions.

\begin{proof}[Proof of Lemma~\ref{lem:Pi E}]
 To begin with, notice that \eqref{eqn:Pi E} is equivalent to: 
\[
\lambda,\mu,\sigma \in F, d(\lambda) = d(\mu), s(\lambda) = s(\mu), 
\text{ and }(\alpha,\beta) \in \Lmin(\mu,\sigma)\text{ imply }\lambda\alpha \in F.
\]

Let $N := \bigvee_{\lambda \in E} d(\lambda)$. Let $E_0 := E$, and let
\[
\begin{split} E_1 := 
\{&\lambda_1(0,d(\lambda_1))\lambda_2(d(\lambda_1),d(\lambda_2))
\cdots\lambda_j(d(\lambda_{j-1}),d(\lambda_j)) :  \lambda_l \in 
\vee E_0, \\
& d(\lambda_l) \le d(\lambda_{l+1}), s(\lambda_l) = 
r(\lambda_{l+1}(d(\lambda_l),d(\lambda_{l+1})))\text{ for } 1 \le l \le j\}.
\end{split}
\]
The set $E_1$ is finite because $\vee E_0$ is finite.  Furthermore $E_1$ contains $E = E_0$ by definition. 
Suppose that $\lambda
\in E_1$. Then $d(\lambda) = d(\lambda_j)$ for some $\lambda_j \in 
\vee E_0$, so
$d(\lambda) \le N$. If $\lambda,\mu,\sigma \in E_0$ with $d(\lambda) 
= d(\mu)$ and
$s(\lambda)=s(\mu)$, and if $(\alpha,\beta) \in \Lmin(\mu,\sigma)$, then
$\lambda,\mu\alpha\in\vee E_0$ and hence $\lambda\alpha \in E_1$.

Iteratively construct sets $E_i \subset \L$, $i \ge 2$ by
\[
\begin{split} E_{i} :=
\{&\lambda_1(0,d(\lambda_1))\cdots 
\lambda_j(d(\lambda_{j-1}),d(\lambda_j)) : \lambda_l \in \vee E_{i-1},\\
&d(\lambda_l) \le d(\lambda_{l+1}), s(\lambda_l) =
r(\lambda_{l+1}(d(\lambda_l),d(\lambda_{l+1})))\text{ for }1 \le l \le j\}.
\end{split}
\] We claim that for all $i \ge 2$,
\begin{itemize}
\item[(a)] $E_i$ is finite;
\item[(b)] $E_{i-1} \subset E_i$;
\item[(c)] $d(\lambda) \le N$ for all $\lambda \in E_i$;
\item[(d)] if $\lambda,\mu,\sigma\in E_{i-1}$ satisfy $d(\lambda) = 
d(\mu)$, $s(\lambda) =
s(\mu)$, and if $(\alpha,\beta) \in \Lmin(\mu,\sigma)$, then 
$\lambda\alpha \in E_i$; and
\item[(e)] If $E_{i-1} \not= E_i$, then $\min_{\lambda \in E_i 
\setminus E_{i-1}}
|d(\lambda)| > \min_{\mu \in E_{i-1} \setminus E_{i-2}} |d(\mu)|$.
\end{itemize} 
Once we have established (a)--(e), conditions (b), (c) and
(e) combine to ensure that $E_{|N|+1} = E_{|N|}$. With $F := 
E_{|N|}$, it then follows that $E \subset
F$ by (b), $F$ is finite by (a), and $F$ satisfies \eqref{eqn:Pi E} by
(d).

Let $h \ge 1$ and suppose that (a)--(d) hold for $i = h$. We will 
show that (a)--(d)
hold for $i = h+1$. Since we have already established (a)--(d) for 
$i=1$, (a)--(d)
will then follow for all $i \ge 1$ by induction. We have $E_{h+1}$ finite 
because $\L$ is
finitely aligned and $E_{h}$ is finite, giving (a). The inclusion 
$E_{h} \subset \vee
E_{h} \subset E_{h+1}$ gives (b). If $\lambda \in E_{h+1}$, then $d(\lambda) =
d(\lambda_j)$ for some $\lambda_j \in \vee E_{h}$, so $d(\lambda) \le N$ by
definition of $\vee E_{h}$, and (c) for $i =h$. Now suppose that
$\lambda,\mu,\sigma$ and $(\alpha,\beta)$ are as in (d) for $i = 
h+1$. Then $\mu\alpha
\in \vee E_{h}$, and $\lambda\alpha = \lambda(0, d(\lambda))
(\mu\alpha)(d(\mu),d(\mu\alpha))\in E_{h+1}$, giving (d) for $i= h+1$.

To establish (e), suppose that $i \ge 2$ and
$\lambda \in E_i \setminus E_{i-1}$. Then
\[
\lambda = 
\lambda_1(0,d(\lambda_1))\cdots\lambda_j(d(\lambda_{j-1}),d(\lambda_j))
\] where each $\lambda_l \in \vee E_{i-1}$. If every $\lambda_l \in 
E_{i-1}$, then each
$\lambda_l$ may be written as
\[
\lambda_l = 
\lambda_{l,1}(0,d(\lambda_{l,1}))\cdots\lambda_{l,h_l}(d(\lambda_{l,h_l-1}),
d(\lambda_{l,h_l}))
\] where each $\lambda_{l,m} \in \vee E_{i-2}$, and then
\[
\lambda =
\lambda_{1,1}(0,d(\lambda_{1,1}))
\lambda_{1,2}(d(\lambda_{1,1}), d(\lambda_{1,2}))
\cdots\lambda_{j,h_j}(d(\lambda_{j,h_j-1}),d(\lambda_{j,h_j}))
\] 
belongs to $E_{i-1}$ contradicting $\lambda \in E_i \setminus 
E_{i-1}$. Hence there
must be some $l$ such that $\lambda_l \in (\vee E_{i-1}) \setminus E_{i-1}$. By
definition of $\vee E_{i-1}$, there exists $G \subset 
E_{i-1}$ such that
$\lambda_l \in \MCE(G)$. Furthermore, $d(\lambda_l) > d(\sigma)$ for 
all $\sigma \in G$,
for if not we have $\lambda_l \in G \subset E_{i-1}$. If $G \subset 
E_{i-2}$, then
$\lambda_l\in E_{i-1}$, so there exists $\sigma \in (G \setminus 
E_{i-2}) \subset
(E_{i-1} \setminus E_{i-2})$. Hence $|d(\lambda)| \ge |d(\lambda_l)| > 
|d(\sigma)| \ge \min_{\mu \in E_{i-1}\setminus E_{i-2}} |d(\mu)|$,
proving the claim.

Now suppose that $F$ is any finite set containing $E$ and satisfying 
\eqref{eqn:Pi E}. Then $M^t_F$ is a finite-dimensional subspace of
$C^*(\L)^\gamma$ which is closed under taking adjoints.
Hence we need only check that $M^t_{\Pi E}$ is closed under multiplication. 
But if $t_\lambda t^*_\mu$ and $t_\sigma t^*_\tau$ are generators of
$M^t_{\Pi E}$, then $\lambda,\mu,\sigma,\tau$ are
as in \eqref{eqn:Pi E}. Since
\[ t_\lambda t^*_\mu t_\sigma t^*_\tau =
\sum_{(\alpha,\beta)\in\Lmin(\mu,\sigma)}t_{\lambda\alpha}t^*_{\tau\beta},
\]
and since each $\lambda\alpha$ and each $\tau\beta$ belong to $F$ by \eqref{eqn:Pi E}, 
it follows that $t_\lambda t^*_\mu t_\sigma t^*_\tau \in M^t_F$.
\end{proof}

The intersection of a family of sets satisfying \eqref{eqn:Pi E}
also satisfies \eqref{eqn:Pi E}, so we can make the following definition.

\begin{dfn}\label{dfn:Pi E} For any $\L$ and $E$, we define $\Pi E$ to 
be the smallest set 
containing $E$ which satisfies \eqref{eqn:Pi E}; that is
\[
\textstyle\Pi E := \bigcap\{F \subset \L : E \subset F\text{ and $F$ satisfies
\eqref{eqn:Pi E}}\}.
\]
\end{dfn}

\begin{rmk}\label{rmk:consequences} The following consequences of 
Lemma~\ref{lem:Pi E}
will prove useful.
\begin{itemize}
\item[(i)] $\Pi E$ is finite.
\item[(ii)] For $\rho,\xi\in\Pi E$ with $d(\rho)=d(\xi)$ and 
$s(\rho)=s(\xi)$, and for all $\nu \in s(\rho)\L$,
\[
\rho\nu\in\Pi E\text{ if and only if }\xi\nu\in\Pi E:
\] the ``if'' direction follows from \eqref{eqn:Pi E} with $\lambda 
=\rho$, $\mu = \xi$,
and $\sigma=\tau = \xi\nu$, and the ``only if'' direction follows 
from \eqref{eqn:Pi E}
with $\lambda=\mu = \rho\nu$, $\sigma = \rho$, and $\tau = \xi$.
\item[(iii)] If $\rho,\xi \in \Pi E$ and $(\alpha,\beta) \in 
\Lmin(\rho,\xi)$, then
\eqref{eqn:Pi E} with $\lambda = \mu = \rho$ and $\sigma = \tau = 
\xi$ gives $\rho\alpha
= \xi\beta \in \Pi E$; that is to say, $\Pi E$ is closed under taking 
minimal common
extensions, so $\Pi E = \vee(\Pi E)$.
\end{itemize}
\end{rmk}

The next step is to find a family of matrix units for $M^t_{\Pi E}$. 
The trick is first to expess each $t_v$ as a sum of orthogonalised range
projections associated to paths in $\Pi E$.

\begin{prop} \label{prp:partition} 
For each $\lambda \in \Pi E$, define
\[ Q(t)^{\Pi E}_{\lambda} := t_\lambda t^*_\lambda 
\prod_{\substack{\lambda\nu \in \Pi E
\\ d(\nu) > 0}} (t_\lambda t_\lambda^* - t_{\lambda\nu}t^*_{\lambda\nu}).
\] Then $\{Q(t)^{\Pi E}_\lambda : \lambda \in \Pi E\}$ is a family of mutually orthogonal
projections such that
\begin{equation}\label{eqn:Pi E sum}
\prod_{\lambda\in v \Pi E}(t_v-t_\lambda t^*_\lambda)+\sum_{\mu 
\in v \Pi E}
Q(t)^{\Pi E}_\mu = t_v
\end{equation} for all $v \in r(\Pi E)$.
\end{prop}
\begin{proof} 
Fix $v\in r(\Pi E)$. Any 
$G\subset \L$
satisfies \eqref{eqn:Pi E} if and only if $G\cup \{v\}$ satisfies 
\eqref{eqn:Pi E}.
Hence, by Definition~\ref{dfn:Pi E}, $(\Pi E)\cup \{v\}=\Pi(E\cup \{v\})$.

If $v\in\Pi E$, then $\prod_{\lambda\in v \Pi E}(t_v-t_\lambda 
t^*_\lambda)=0$, so
setting $F:=v \Pi E$, the left-hand side of \eqref{eqn:Pi E sum} is equal to
$\sum_{\lambda\in F}Q(t)^F_{\lambda}$.

On the other hand, if $v\not\in \Pi E$, then with $F:=v((\Pi E)\cup\{ 
v\})$, we have
\[ Q(t)^F_\lambda = Q(t)^{(\Pi E)\cup \{v\}}_\lambda=Q(t)^{\Pi E}_\lambda
\] for all $\lambda\in v(\Pi E)$. Furthermore,
\[ Q(t)^F_v=\prod_{\lambda\in v \Pi E}(t_v-t_\lambda t^*_\lambda).
\] So the left-hand side of \eqref{eqn:Pi E sum} is once again equal to
$\sum_{\lambda\in F}Q(t)^F_{\lambda}$.

In either case, $F=\vee F$ and $\lambda\in F\implies r(\lambda)\in 
F$. Under the
identification of finitely aligned product systems of graphs over 
$\NN^k$ with finitely
aligned $k$-graphs (see \cite[Example~3.5]{RS1}), the proof of
\cite[Proposition~8.6]{RS1}  with its first sentence removed now 
proves our result.
\end{proof}

\begin{rmk}\label{rmk:empty products}
For $\lambda \in \Pi E$, we have
\begin{equation}\label{eqn:the trick} \begin{split}
Q(t)^{\Pi E}_\lambda 
&= t_\lambda t^*_\lambda
\prod_{\substack{\lambda\nu \in \Pi E \\ d(\nu) > 0}}( t_\lambda( 
t_{s(\lambda)} - t_\nu
t^*_\nu) t^*_\lambda) \\
&= t_\lambda \Big(\prod_{\substack{\lambda\nu 
\in \Pi E \\ d(\nu) >
0}} (t_{s(\lambda)} - t_\nu t^*_\nu) \Big)t^*_\lambda
\end{split}\end{equation} 
because $t^*_\lambda t_\lambda = t_{s(\lambda)}$.
\end{rmk}

\begin{cor}\label{cor:Qs add to mu} Let $\mu\in \Pi E$. Then $t_\mu t^*_\mu
=\sum_{\mu\nu\in\Pi E}Q(t)^{\Pi E}_{\mu\nu}$.
\end{cor}

\begin{proof} First notice that
\[ t_\mu t^*_\mu = t_\mu t^*_\mu t_{r(\mu)} = t_\mu t^*_\mu 
\Big(\prod_{\lambda\in
r(\mu)\Pi E}(t_{r(\mu)}-t_\lambda t^*_\lambda)+\sum_{\sigma \in 
r(\mu) \Pi E} Q(t)^{\Pi E}_\sigma\Big)
\] by Proposition~\ref{prp:partition}. By definition of $Q(t)^{\Pi 
E}_{\mu\nu}$, we have
$t_\mu t^*_\mu \ge Q(t)^{\Pi E}_{\mu\nu}$ for all $\nu$, so 
it suffices to show
that
\begin{itemize}
\item[(i)]
$t_\mu t^*_\mu \prod_{\lambda\in r(\mu)\Pi E}(t_{r(\mu)}-t_\lambda 
t^*_\lambda) = 0$; and
\item[(ii)] for $\sigma \in \Pi E$ with $\sigma(0,d(\mu)) \not= \mu$, 
we have $t_\mu
t^*_\mu Q(t)^{\Pi E}_\sigma = 0$.
\end{itemize} Claim (i) is straightforward because $\mu \in r(\mu) 
\Pi E$, and hence
\[ t_\mu t^*_\mu \prod_{\lambda\in r(\mu)\Pi E}(t_{r(\mu)}-t_\lambda 
t^*_\lambda) \le
t_\mu t^*_\mu (t_{r(\mu)} - t_\mu t^*_\mu) = 0.
\]

It remains to prove Claim (ii). But for $\sigma$ as in Claim (ii), 
$(\alpha,\beta) \in
\Lmin(\mu,\sigma)$ implies $d(\beta) > 0$, and the definition of $\Pi 
E$ ensures that
$\sigma\beta \in \Pi E$. Hence
\begin{align*} t_\mu t^*_\mu &Q(t)^{\Pi E}_\sigma \\ &= t_\mu t^*_\mu t_\sigma
t^*_\sigma \prod_{\substack{\sigma\nu \in \Pi E \\ d(\nu) > 0}} 
(t_\sigma t^*_\sigma -
t_{\sigma\nu}t^*_{\sigma\nu}) \\ &= \Big(\sum_{(\alpha,\beta) \in \Lmin(\mu,\sigma)}
t_{\sigma\beta} t^*_{\sigma\beta}\Big)\Big(\prod_{\substack{\sigma\nu 
\in \Pi E \\
d(\nu) > 0}} (t_\sigma t^*_\sigma - t_{\sigma\nu}t^*_{\sigma\nu})\Big) \\ &=
\sum_{(\alpha,\beta)\in 
\Lmin(\mu,\sigma)}\Big(t_{\sigma\beta}t^*_{\sigma\beta}(t_\sigma
t^*_\sigma - t_{\sigma\beta} t^*_{\sigma\beta}) 
\prod_{\substack{\sigma\nu \in \Pi E \setminus\{\sigma\beta\} \\
d(\nu) > 0}} (t_\sigma t^*_\sigma - 
t_{\sigma\nu}t^*_{\sigma\nu})
\Big) \\ &= 0,
\end{align*} establishing Claim (ii).
\end{proof}

\begin{dfn}\label{dfn:thetas} For $\lambda,\mu \in \Pi E$ with 
$d(\lambda) = d(\mu)$ and
$s(\lambda) = s(\mu)$, define
$\MU(t,E;\lambda,\mu) := Q(t)^{\Pi E}_\lambda t_\lambda t^*_\mu$.
\end{dfn}

\begin{prop}\label{prp:matrix units} The set
\[
\{\MU(t,E;\lambda,\mu) : \lambda,\mu \in \Pi E, d(\lambda) = d(\mu), 
s(\lambda) =
s(\mu)\}
\] is a collection of partial isometries which span $M^t_{\Pi E}$ and satisfy
\begin{itemize}
\item[(i)] $\big(\MU(t,E;\lambda,\mu)\big)^* =  \MU(t,E;\mu,\lambda)$; and
\item[(ii)] $\MU(t,E;\lambda,\mu) \MU(t,E;\sigma,\tau) = \delta_{\mu,\sigma} \MU(t,E;\lambda,\tau)$.
\end{itemize}
\end{prop}

To prove Proposition~\ref{prp:matrix units} we need to establish two lemmas.

\begin{lemma}\label{lem:claim1} Let $\lambda,\mu \in \Pi E$ with 
$d(\lambda) = d(\mu)$
and $s(\lambda) = s(\mu)$. Then
\[
\MU(t,E;\lambda,\mu) = t_\lambda\Big(\prod_{\substack{\lambda\nu \in 
\Pi E\\ d(\nu) >
0}}(t_{s(\lambda)} - t_\nu t^*_\nu)\Big)t^*_\mu = t_\lambda t^*_\mu 
Q(t)^{\Pi E}_\mu.
\]
\end{lemma}
\begin{proof} We begin by calculating:
\begin{align}
\MU(t,E;\lambda,\mu) &= Q(t)^{\Pi E}_\lambda t_\lambda t^*_\mu \notag\\ &=
t_\lambda\Big(\prod_{\substack{\lambda\nu \in \Pi E\\ d(\nu) > 
0}}(t_{s(\lambda)} -
t_\nu t^*_\nu)\Big) t^*_\lambda t_\lambda t^*_\mu \quad\text{by \eqref{eqn:the
trick}}\notag\\ &= t_\lambda\Big(\prod_{\substack{\lambda\nu \in \Pi 
E\\ d(\nu) >
0}}(t_{s(\lambda)} - t_\nu t^*_\nu)\Big)t^*_\mu, \label{eqn:first =}
\end{align} which establishes the first equality. For the second 
equality, we continue
the calculation as follows:
\begin{flalign*} &&\MU(t,E;\lambda,\mu) &= 
t_\lambda\Big(\prod_{\substack{\lambda\nu \in
\Pi E\\ d(\nu) > 0}}(t_{s(\lambda)} - t_\nu t^*_\nu)\Big)t^*_\mu
\quad\text{by \eqref{eqn:first =}}&\\ &&&= 
t_\lambda\Big(\prod_{\substack{\mu\nu \in \Pi
E\\ d(\nu) > 0}}(t_{s(\lambda)} - t_\nu t^*_\nu)\Big)t^*_\mu
\quad\text{by Remark~\ref{rmk:consequences}(ii)} &\\ &&&= t_\lambda 
t^*_\mu \Big( t_\mu
\prod_{\substack{\mu\nu \in \Pi E\\ d(\nu) > 0}}(t_{s(\lambda)} - 
t_\nu t^*_\nu)t^*_\mu
\Big) &\\ &&&= t_\lambda t^*_\mu Q(t)^{\Pi E}_\mu \quad\text{by 
\eqref{eqn:the trick}.}&\qed
\end{flalign*}
\renewcommand{\qed}{}\end{proof}

\begin{lemma}\label{lem:claim3} Let $\lambda,\mu \in \Pi E$ with 
$d(\lambda) = d(\mu)$
and $s(\lambda) = s(\mu)$. Then
\[ t_\lambda t^*_\mu = \sum_{\lambda\nu \in \Pi E} \MU(t,E;\lambda\nu, \mu\nu).
\]
\end{lemma}
\begin{proof} Just calculate
\begin{flalign*} &&t_\lambda t^*_\mu &= t_\lambda t^*_\mu t_\mu t^*_\mu &\\
&&&=t_\lambda t^*_\mu \Big(\sum_{\mu\nu \in \Pi E} Q(t)^{\Pi E}_{\mu\nu}\Big)
\quad\text{by Corollary~\ref{cor:Qs add to mu}} &\\ &&&=\sum_{\mu\nu \in \Pi E}
\Big(t_\lambda t^*_\mu t_{\mu\nu}\Big(\prod_{\substack{\mu\nu\nu' \in 
\Pi E \\ d(\nu') >
0}}( t_{s(\nu)} - t_{\nu'} t^*_{\nu'}) t^*_{\mu\nu}\Big)\Big) \quad\text{by
\eqref{eqn:the trick}} &\\ &&&= \sum_{\lambda\nu \in \Pi E}
\Big(t_{\lambda\nu}\Big(\prod_{\substack{\lambda\nu\nu' \in \Pi E \\ 
d(\nu') > 0}}(
t_{s(\nu)} - t_{\nu'} t^*_{\nu'})\Big) t^*_{\mu\nu}\Big) &\\ 
&&&\phantom{=}\qquad\text{by two
applications of Remark~\ref{rmk:consequences}(ii)}&\\ &&&= 
\sum_{\lambda\nu \in \Pi E}
\MU(t,E;\lambda\nu, \mu\nu) \quad\text{by Lemma~\ref{lem:claim1}} &\qed
\end{flalign*}
\renewcommand\qed{}\end{proof}

\begin{proof}[Proof of Proposition~\ref{prp:matrix units}]
The $\MU(t,E;\lambda,\mu)$ are clearly partial isometries. It follows from
Lemma~\ref{lem:claim3} that they span $M^t_{\Pi E}$. It remains to 
show that the $\MU(t,E;\lambda,\mu)$ satisfy  (i) and (ii).

Let $\lambda,\mu \in \Pi E$ with $d(\lambda) = d(\mu)$ 
and $s(\lambda) = s(\mu)$. Since the $Q(t)^{\Pi E}_\lambda$ are projections by
Proposition~\ref{prp:partition}, we can and use 
Lemma~\ref{lem:claim1} to calculate
\[ 
(\MU(t,E;\lambda,\mu))^* = (Q(t)^{\Pi E}_\lambda t_\lambda 
t^*_\mu)^* = t_\mu
t^*_\lambda Q(t)^{\Pi E}_\lambda = \MU(t,E;\mu,\lambda).
\] 
Furthermore, if $\sigma,\tau$ also belong to $\Pi E$ with 
$d(\sigma) = d(\tau)$ and
$s(\sigma) = s(\tau)$, then
\begin{flalign*} &&\MU(t,E;\lambda,\mu) \MU(t,E;\sigma,\tau) &= 
t_\lambda t^*_\mu
Q(t)^{\Pi E}_\mu Q(t)^{\Pi E}_\sigma t_\sigma t^*_\tau \quad\text{by
Lemma~\ref{lem:claim1}} &\\ &&&=\delta_{\mu,\sigma} t_\lambda t^*_\mu 
Q(t)^{\Pi E}_\mu
t_\mu t^*_\tau \quad\text{by Proposition~\ref{prp:partition}} &\\ &&&=
\delta_{\mu,\sigma} Q(t)^{\Pi E}_\lambda t_\lambda t^*_\mu t_\mu 
t^*_\tau \quad\text{by
Lemma~\ref{lem:claim1}} &\\ &&&= \delta_{\mu,\sigma} Q(t)^{\Pi 
E}_\lambda t_\lambda
t^*_\tau \quad\text{since $s(\lambda) = s(\mu)$} &\\ &&&=\delta_{\mu,\sigma}
\MU(t,E;\lambda,\tau) &\qed
\end{flalign*}
\renewcommand\qed{}\end{proof}

We now need to say which pairs $\lambda,\mu$ satisfy $\MU(t,E;\lambda,\mu) \not= 0$. 

\begin{notation} \label{not:T(lambda)}
For $\lambda,\mu \in \Pi E$ with $s(\lambda) = s(\mu) = v$ and
$d(\lambda) = d(\mu) = n$, Remark~\ref{rmk:consequences}(ii) ensures that
\[
\{\nu \in v\L : d(\nu) > 0, \lambda\nu \in \Pi E\} 
= \{\nu \in v\L : d(\nu) > 0, \mu\nu \in \Pi E\}.
\] 
We denote this set by $T^{\Pi E}(n, v)$. For 
convenience, for $\lambda \in \Pi E$, we write
$T(\lambda)$ for $T^{\Pi E}(d(\lambda), s(\lambda))$.
\end{notation}

\begin{prop}\label{prp:nonzero matrix units}
Suppose that $t_v \not= 0$ for all $v \in \L^0$. 
Then 
\[
\MU(t,E;\lambda,\mu) = 0\text{ if and only if 
$T(\lambda)$ is exhaustive.}
\]
\end{prop}

To prove Proposition~\ref{prp:nonzero matrix units}, we need a definition and two 
lemmas.

\begin{dfn}\label{dfn:xi extendors} For each $n \in \NN^k$ and $v \in 
\L^0$ with $T^{\Pi E}(n,v)$ non-exhaustive, fix $\xi^{\Pi E}(n,v) \in v\L$ such that 
$\Lmin(\xi^{\Pi E}(n,v),\nu) = \emptyset$ for all $\nu \in T^{\Pi E}(n,v)$. Again for
convenience, we will write $\xi_\lambda$ in place of $\xi^{\Pi
E}(d(\lambda),s(\lambda))$ for $\lambda \in \Pi E$.
\end{dfn}

\begin{lemma}\label{lem:claim4} For each $\lambda \in \Pi E$ such 
that $T(\lambda)$ is
not exhaustive, $t_{\lambda\xi_\lambda} t^*_{\lambda\xi_\lambda} 
\le Q(t)^{\Pi
E}_\lambda$.
\end{lemma}
\begin{proof} Set $\xi = \xi_\lambda$, and calculate
\begin{flalign*} &&t_{\lambda\xi} t^*_{\lambda\xi} Q(t)^{\Pi E}_\lambda &=
t_{\lambda\xi} t^*_{\lambda\xi} t_\lambda t^*_\lambda 
\prod_{\substack{\lambda\nu \in
\Pi E \\ d(\nu) > 0}}(t_\lambda t^*_\lambda - t_{\lambda\nu} 
t^*_{\lambda\nu}) &\\ &&&=
\prod_{\substack{\lambda\nu \in \Pi E \\ d(\nu) > 0}}\big(t_{\lambda\xi}
t^*_{\lambda\xi}(t_\lambda t^*_\lambda - t_{\lambda\nu} 
t^*_{\lambda\nu})\big) &\\ &&&=
\prod_{\substack{\lambda\nu \in \Pi E \\ d(\nu) > 0}}\Big(t_{\lambda\xi}
t^*_{\lambda\xi} - \sum_{(\alpha,\beta) \in \Lmin(\lambda\xi, \lambda\nu)}
t_{\lambda\nu\beta} t^*_{\lambda\nu\beta}\Big) &\\ &&&= 
\prod_{\substack{\lambda\nu \in
\Pi E \\ d(\nu) > 0}} t_{\lambda\xi} t^*_{\lambda\xi} &\\ 
&&&\phantom{=}\qquad\text{since each 
$\Lmin(\lambda\xi,\lambda\nu) = \Lmin(\xi,\nu) = \emptyset$}&\\ 
&&&\phantom{=}\qquad\text{by choice of $\xi = \xi_\lambda$} &\\ 
&&&= t_{\lambda\xi} t^*_{\lambda\xi}&\qed
\end{flalign*}
\renewcommand\qed{}\end{proof}

\begin{lemma}\label{lem:claim5} Let $\lambda\in \Pi E$ and suppose 
that $T(\lambda)$ is
not exhaustive. Let $\sigma,\tau \in \Pi E$ with $d(\sigma) = 
d(\tau)$ and $s(\sigma) =
s(\tau)$. Then
\[ t_{\lambda\xi_\lambda} t^*_{\lambda\xi_\lambda}\,\MU(t,E;\sigma,\tau) =
\delta_{\lambda,\sigma} t_{\lambda\xi_\lambda} t^*_{\tau\xi_\lambda}.
\]
\end{lemma}
\begin{proof} Set $\xi = \xi_\lambda$ and calculate
\begin{flalign*} 
&& t_{\lambda\xi} t^*_{\lambda\xi}\,\MU(t,E;\sigma,\tau) &=
t_{\lambda\xi} t^*_{\lambda\xi}\,Q(t)^{\Pi E}_{\sigma} t_\sigma 
t^*_\tau &\\
&&&= t_{\lambda\xi} t^*_{\lambda\xi}Q(t)^{\Pi E}_\lambda Q(t)^{\Pi E}_{\sigma} t_\sigma 
t^*_\tau \qquad\text{by Lemma~\ref{lem:claim4}}&\\
 &&&= \delta_{\lambda,\sigma} t_{\lambda\xi}
t^*_{\lambda\xi}\, Q(t)^{\Pi E}_\lambda t_\lambda t^*_\tau
\qquad\text{by Proposition~\ref{prp:partition}} &\\ 
&&&= \delta_{\lambda,\sigma} t_{\lambda\xi}
t^*_{\tau\xi} 
\qquad\text{by Lemma~\ref{lem:claim4}} &\qed
\end{flalign*}
\renewcommand\qed{}\end{proof}

\begin{proof}[Proof of Proposition~\ref{prp:nonzero matrix units}] 
For the ``if'' direction, note that $T(\lambda)$ is 
certainly finite and if it is also exhaustive then
\[
\MU(t,E;\lambda,\mu) = t_\lambda \Big(\prod_{\nu \in T(\lambda)} 
(t_{s(\lambda)} - t_\nu t^*_\nu)\Big)t_\mu = 0
\]
by Definition~\ref{dfn:CKfamily}(iv). 
For the ``only if'' direction, suppose that $\lambda,\mu \in \Pi E$ with 
$d(\lambda) = d(\mu)$
and $s(\lambda) = s(\mu)$, and suppose that $T(\lambda)$ is not 
exhaustive. Then
Lemma~\ref{lem:claim5} ensures that
\[ t_{\lambda\xi_\lambda} t^*_{\lambda\xi_\lambda} \MU(t,E;\lambda,\mu) =
t_{\lambda\xi_\lambda} t^*_{\lambda\xi_\lambda},
\] which is nonzero because each $t_v \not=0$. 
Hence $\MU(t,E;\lambda,\mu) \not= 0$.
\end{proof}

\begin{cor}\label{cor:nonzero MUs} 
Suppose that $t_v \not= 0$ for all $v \in \L^0$. Suppose $\lambda,\mu \in \Pi E$ with $d(\lambda) = d(\mu)$ and
$s(\lambda) = s(\mu)$. Then $\MU(t,E;\lambda,\mu) = 0$ if and only if
$\MU(s,E;\lambda,\mu) = 0$.
\end{cor}
\begin{proof} We know from the boundary path representation that each 
$s_v$ is nonzero. The result then follows from 
Proposition~\ref{prp:nonzero matrix units} applied to both
$\{s_\lambda\}$ and $\{t_\lambda\}$.
\end{proof}

\begin{proof}[Proof of Theorem~\ref{thm:faithful on core}]
 Since 
\[
C^*(\L)^\gamma = \clsp\{s_\lambda s^*_\mu : 
\lambda,\mu \in \L, d(\lambda) = d(\mu)\},
\]
we have
\[ 
C^*(\L)^\gamma = \overline{\bigcup_{E \subset \L \text{ finite}} 
M^s_{\Pi E}}.
\] 
Since each $M^s_{\Pi E}$ is finite-dimensional, 
it follows that $C^*(\L)^\gamma$ is AF. Furthermore, since
$\pi_t(\MU(s,E;\lambda,\mu)) = \MU(t,E;\lambda,\mu)$ for all finite 
$E \subset \L$ and $\MU(t,E;\lambda,\mu) \in M^t_{\Pi E}$,
Corollary~\ref{cor:nonzero MUs} ensures that $\pi_t$ maps nonzero matrix units
$\MU(s,E;\lambda,\mu)$ to nonzero matrix units 
$\MU(t,E;\lambda,\mu)$, and hence is
faithful on each $M^s_{\Pi E}$. The result now follows from 
\cite[Lemma~1.3]{ALNR}.
\end{proof}

\section{The uniqueness theorems} \label{sec:theorems}

Write $\Phi$ for the linear map from $C^*(\L)$ to $C^*(\L)^\gamma$ 
obtained by averaging
over the gauge action; that is, $\Phi(a) := \int_{\TT^k} \gamma_z(a) 
dz$. The map $\Phi$
is faithful on positive elements and satisfies $\Phi(s_\lambda s^*_\mu) =
\delta_{d(\lambda),d(\mu)} s_\lambda s^*_\mu$.

\begin{prop}\label{prp:standard argument} Let $(\L,d)$ be a finitely 
aligned $k$-graph.
Suppose that $\pi$ is a homomorphism of $C^*(\L)$ such that 
$\pi(s_v) \not= 0$ for
all $v \in \L^0$ and
\begin{equation}\label{eqn:pi circ Phi norm reducing}
\|\pi(\Phi(a))\| \le \|\pi(a)\|\text{ for all $a \in C^*(\L)$.}
\end{equation} Then $\pi$ is injective.
\end{prop}

\begin{proof} Equation~\eqref{eqn:pi circ Phi norm reducing},
Theorem~\ref{thm:faithful on core}, and the properties of $\Phi$ 
show that $\pi(a^*a)
= 0 \implies a^*a = 0$.
\end{proof}

\subsection{The gauge-invariant uniqueness theorem}

\begin{theorem}\label{thm:giut} Let $(\L,d)$ be a finitely aligned 
$k$-graph, and let
$\pi$ be a homomorphism of $C^*(\L)$. Suppose that there is a 
strongly continuous
action $\theta : \TT^k \to \Aut\big(C^*(\{\pi(s_\lambda) : \lambda 
\in \L\})\big)$ such
that $\theta_z \circ \pi = \pi \circ \gamma_z$ for all $z \in \TT^k$. 
If $\pi(s_v)\neq
0$ for all $v\in\L^0$, then $\pi$ is injective.
\end{theorem}

\begin{proof} Averaging over $\theta$ is norm-decreasing and 
implements $\pi(a) \mapsto
\pi(\Phi(a))$. Hence Equation~\eqref{eqn:pi circ Phi norm reducing} 
holds, and the
result follows from Proposition~\ref{prp:standard argument}.
\end{proof}

\begin{cor}[The gauge-invariant uniqueness theorem] \label{cor:giut} 
Let $(\L,d)$ be a
finitely aligned $k$-graph. There exists a Cuntz-Krieger $\L$-family 
$\{t_\lambda :
\lambda \in \L\}$ such that $t_v \not= 0$ for every $v \in \L^0$, and 
such that there
exists a strongly continuous action $\theta : \TT^k \to 
\Aut(C^*(\{t_\lambda : \lambda
\in \Lambda\}))$ satisfying $\theta_z(t_\lambda) = z^{d(\lambda)} 
t_\lambda$ for all
$\lambda \in \Lambda$. Furthermore, any two such families generate canonically
isomorphic $C^*$-algebras.
\end{cor}

\begin{proof} Proposition~\ref{prp:bdry repn nonzero} shows that there is a
Cuntz-Krieger $\L$-family consisting of nonzero partial isometries. 
It follows that each
$s_v \in C^*(\L)$ is nonzero, so $t_\lambda := 
s_\lambda$ and $\theta :=
\gamma$ gives existence. The last statement of the corollary follows from
Theorem~\ref{thm:giut}.
\end{proof}

Recall from \cite{KP} that if $(\L_1,d_1)$ is a $k_1$-graph
and $(\L_2,d_2)$ is a $k_2$-graph,
then the pair $(\L_1 \times \L_2, d_1 \times d_2)$ is a $(k_1 + 
k_2)$-graph. It is easy
to check that if $\L_1$ and $\L_2$ are finitely aligned, then so is 
$\L_1 \times \L_2$.

\begin{cor}\label{cor:cartesian products} Let $\L_1$ be a finitely aligned
$k_1$-graph and let $\L_2$ be a finitely
aligned $k_2$-graph. Then $C^*(\L_1 \times \L_2)$ is
canonically isomorphic to $C^*(\L_1) \otimes C^*(\L_2)$.
\end{cor}

\begin{proof} Implicit in the statement of the corollary is that all 
tensor products of
$C^*(\L_1)$ and $C^*(\L_2)$ coincide. The 
bilinearity of tensor
products ensures that $\{s_{\lambda_1} \otimes s_{\lambda_2} : 
(\lambda_1, \lambda_2)
\in \L_1 \times \L_2\}$ is a Cuntz-Krieger $(\L_1 \times 
\L_2)$-family regardless of the
tensor product in question. Seperate arguments for the spatial tensor 
product and the
universal tensor product show that for either one, the formula
\[
\theta_{z}(s_{\lambda_1} \otimes s_{\lambda_2}) := 
\big(z_1^{d(\lambda_1)_1} \cdots
z_{k_1}^{d(\lambda_1)_{k_1}} z_{k_1 + 1}^{d(\lambda_2)_1} \cdots z_{k_1 +
k_2}^{d(\lambda_2)_{k_2}}\big) s_{\lambda_1} \otimes s_{\lambda_2}
\] 
extends to a strongly continuous action $\theta$ of $\TT^{k_1 + k_2}$ on
$C^*(\{s_{\lambda_1} \otimes s_{\lambda_2} : (\lambda_1, \lambda_2) 
\in \L_1 \times
\L_2\})$ which is equivariant with the gauge action on $C^*(\L_1 
\times \L_2)$. The vertex projections $s_{v_1} \otimes s_{v_2}$ are all nonzero
because each $s_{v_1}$ is nonzero and each $s_{v_2}$ is nonzero. 
Corollary~\ref{cor:giut} shows that the two tensor products coincide, and
Theorem~\ref{thm:giut} shows they are canonically isomorphic to $C^*(\L_1 \times \L_2)$.
\end{proof}

\subsection{The Cuntz-Krieger uniqueness theorem}

\begin{theorem}\label{thm:CKut} Let $(\L,d)$ be a finitely aligned 
$k$-graph, and
suppose that
\begin{equation}\label{eqn:aperiodicity}
\begin{split}
\text{for each $v \in \L^0$ there}&\text{ exists $x \in v 
\uL(\infty)$ such that} \\
&\text{$\lambda,\mu \in \L v$ and $\lambda \not= \mu$ imply $\lambda 
x \not= \mu x$.}
\end{split}\tag{B}
\end{equation} Suppose that $\pi$ is a homomorphism 
of $C^*(\L)$ such that $\pi(s_v)
\not= 0$ for all $v \in \L^0$. Then $\pi$ is injective.
\end{theorem}

\begin{cor}[The Cuntz-Krieger uniqueness theorem] Let $(\L,d)$ be a 
finitely aligned
$k$-graph which satisfies \eqref{eqn:aperiodicity}. There exists a 
Cuntz-Krieger
$\L$-family $\{t_\lambda : \lambda \in \L\}$ such that $t_v \not= 0$ 
for all $v \in
\L^0$. Furthermore, any two such families generate canonically 
isomorphic $C^*$-algebras.
\end{cor}

\begin{proof} The existence of a nonzero Cuntz-Krieger $\L$-family follows from
Proposition~\ref{prp:bdry repn nonzero}. The last statement of the 
corollary follows
from Theorem~\ref{thm:CKut}.
\end{proof}

The rest of this section is devoted to proving 
Theorem~\ref{thm:CKut}. For the remainder
of this section, let $(\L,d)$ and $\pi$ be as in 
Theorem~\ref{thm:CKut}
and fix a finite set $E \subset \L$ and a linear combination $a = 
\sum_{\lambda,\mu \in
E} a_{\lambda,\mu} s_\lambda s^*_\mu \in C^*(\L)$. Notice that $\Phi(a) =
\sum_{\lambda,\mu \in E, d(\lambda) = d(\mu)} a_{\lambda,\mu} 
s_\lambda s^*_\mu$. Since
$a$ is arbitrary in a dense subset of $C^*(\L)$, if we show that
\[
\|\pi(\Phi(a))\| \le \|\pi(a)\|,
\] then Theorem~\ref{thm:CKut} will follow from 
Proposition~\ref{prp:standard argument}.

For $n \in \NN^k$, define $\Ff_n$ to be the $C^*$-subalgebra of 
$C^*(\L)^\gamma$,
\begin{align*}
\Ff_n 
&:= \clsp\{s_\lambda s^*_\mu : \lambda,\mu \in \uL(n), 
d(\lambda) = d(\mu)\} \\
&\cong \bigoplus_{v \in \L^0,\,m \le n} \Kk(\ell^2(v\uL(n) \cap \L^m))
\end{align*}
where the isomorphism follows from Lemma~\ref{lem:uptoLambda orth range projections}(ii).

\begin{prop}\label{prp:into Ff} There exists ${N_E} \in \NN^k$ and a 
projection $P_E$ such that  $b \mapsto P_E b$ 
is an isomorphism of
$M^s_{\Pi E}$ into $\Ff_{{N_E}}$.
\end{prop}
\begin{proof} Recalling Notation~\ref{not:T(lambda)} and Definition~\ref{dfn:xi extendors}, 
let 
\[\textstyle
{N_E} := \bigvee\{d(\lambda\xi_\lambda)  : \lambda 
\in \Pi E, T(\lambda)\text{ non-exhaustive}\}.
\]
Whenever $T^{\Pi E}(n,v)$ is non-exhaustive, 
$d(\xi^{\Pi E}(n,v)) \le N_E - n$, so let  
let $\nu^{\Pi E}(n,v) \in \uL({{N_E} -n})$ be an extension of $\xi^{\Pi E}(n,v)$. 
That is, for $\lambda \in \Pi E$,  $\nu_\lambda := 
\nu^{\Pi E}(d(\lambda),s(\lambda))$ belongs to
$\uL({{N_E} - d(\lambda)})$ and $\nu_\lambda(0,d(\xi_\lambda) = 
\xi_\lambda$.

Let
\[ 
P_E := \sum_{\substack{\lambda \in \Pi E \\ T(\lambda)\text{ non-exh.}}}
s_{\lambda\nu_\lambda} s^*_{\lambda\nu_\lambda}.
\] 
For all $\lambda \in \Pi E$ with $T(\lambda)$ non-exhaustive,
\[ s_{\lambda\nu_\lambda} s^*_{\lambda\nu_\lambda} \le 
s_{\lambda\xi_\lambda}
s^*_{\lambda\xi_\lambda} \le Q(s)^{\Pi E}_{\lambda}
\] 
by Lemma~\ref{lem:claim5}. Since all the $Q(t)^{\Pi E}_{\lambda}$ 
are mutually orthogonal by Proposition~\ref{prp:partition}, 
it follows that the $s_{\lambda\xi_\lambda} 
s^*_{\lambda\xi_\lambda}$
are mutually orthogonal, as are the $s_{\lambda\nu_\lambda} 
s^*_{\lambda\nu_\lambda}$.
Hence, for all $\lambda \in \Pi E$ with $T(\lambda)$ non-exhaustive,
\begin{equation} \label{eqn:PsubE bigger} 
P_E s_{\lambda\xi_\lambda}
s^*_{\lambda\xi_\lambda} = s_{\lambda\nu_\lambda} 
s^*_{\lambda\nu_\lambda}.
\end{equation} 
If $\lambda,\mu \in \Pi E$ with $d(\lambda) = 
d(\mu)$, $s(\lambda) =
s(\mu)$ and $T(\lambda)$ non-exhaustive, then
\begin{align} 
P_E\,\MU(s,E;\lambda,\mu) & = 
P_E \Big(\sum_{\substack{\sigma \in \Pi E \\ T(\sigma)\text{ non-exh.}}}
s_{\sigma\xi_\sigma} s^*_{\sigma\xi_\sigma}\Big) 
\MU(s,E;\lambda,\mu)  \quad\text{by
\eqref{eqn:PsubE bigger}} \notag \\ &= P_E s_{\lambda\xi_\lambda} 
s^*_{\mu\xi_\lambda}
\quad\text{by Lemma~\ref{lem:claim5}} \notag \\ &= s_{\lambda\nu_\lambda}
s^*_{\mu\nu_\lambda} \quad\text{by \eqref{eqn:PsubE bigger}.}
\label{eqn:PsubE mpctn}
\end{align} 
Lemma~3.6 of \cite{RSY1} says that if $\lambda \in \uL(n)$ and $\mu \in \uL(m)$ then 
$\lambda\mu \in \uL(n+m)$. Hence for all $\lambda \in \Pi E$ such that $T(\lambda)$ is non-exhaustive,
$\lambda\nu_\lambda \in \uL({{N_E}})$. It follows from Proposition~\ref{prp:nonzero matrix units} that $b \mapsto 
P_E b$ sends nonzero
matrix units in $M^s_{\Pi E}$ to nonzero matrix units in 
$\Ff_{{N_E}}$, proving that $b
\mapsto P_E b$ is an isomorphism.
\end{proof}

For $v \in s(\{\nu_\lambda : \lambda \in \Pi E, T(\lambda)\text{ non-exhaustive}\})$, define 
\[
P_v := \sum_{\substack{\lambda\in \Pi E,\, T(\lambda)\text{ non-exh.} \\
 s(\nu_\lambda) = v}} 
s_{\lambda\nu_\lambda} 
s^*_{\lambda\nu_\lambda}, 
\]
so $P_E =
\sum_{v \in s(\{\nu_\lambda : \lambda\in \Pi E,\text{ $T(\lambda)$ non-exh.}\})} P_v$.
In particular $P_v = P_v P_E $, so Equation~\eqref{eqn:PsubE mpctn} gives
\[
P_v \MU(s,E;\lambda,\mu) 
= P_v P_E \MU(s,E;\lambda,\mu)
= P_v s_{\lambda \nu_\lambda} s^*_{\mu \nu_\lambda}
= \delta_{v, s(\nu_\lambda)} s_{\lambda \nu_\lambda} s^*_{\mu \nu_\lambda}.
\]
for all $\lambda,\mu \in \Pi E$ with
$d(\lambda) = d(\mu)$, $s(\lambda) = s(\mu)$, and $T(\lambda) = T(\mu)$ non-exhaustive.
Hence
\[\begin{split}
\MU(s, E;\lambda,\mu) P_v = (P_v \MU(s, E;\mu,\lambda))^* 
&= (\delta_{v,s(\nu_\mu)} s_{\mu \nu_\mu} s^*_{\lambda\nu_\mu})^* \\
&= \delta_{v,s(\nu_\lambda)} s_{\lambda \nu_\lambda} s^*_{\mu\nu_\lambda}\
= P_v \MU(s, E;\lambda,\mu),
\end{split}\]
so each $P_v$ is in the commutant of $M^s_{\Pi E}$.
It follows that there exists a vertex $v_0$
such that
\begin{equation}\label{eqn:norm equality}
\|P_{v_0} \Phi(a)\| = \|P_E \Phi(a)\| = \|\Phi(a)\|
\end{equation} where the second equality follows from 
Proposition~\ref{prp:into Ff}.

\begin{lemma}\label{lem:extensions} Let $\lambda,\mu \in \Pi E$, suppose that
$T(\lambda)$ is not exhaustive, and suppose that $\lambda \not\in 
\mu\L$. Then $\Lmin(\lambda\nu_\lambda, \mu) = \emptyset$.
\end{lemma}
\begin{proof} Suppose for contradiction that $(\eta,\zeta) \in
\Lmin(\lambda\nu_\lambda, \mu)$. Then $\eta = s(\nu_\lambda)$ and 
$\lambda\nu_\lambda
= \mu\zeta$ because $\lambda\nu_\lambda \in \uL({{N_E}})$ and ${N_E} 
\ge d(\mu)$ by
definition. But then with
\[
\alpha := \nu_\lambda(0, (d(\lambda) \vee d(\mu)) - d(\lambda))
\quad\text{and}\quad
\beta := \zeta(0, (d(\lambda) \vee d(\mu)) - d(\mu)),
\] we have $(\alpha,\beta) \in \Lmin(\lambda,\mu)$, and 
$\lambda \not= \mu\mu'$, so $d(\alpha) > 0$; hence
$\alpha \in T(\lambda)$.
Furthermore, $\Lmin(\alpha,\nu_\lambda) \not= \emptyset$  by definition
of $\alpha$, and hence $\Lmin(\xi_\lambda, \alpha) \not= \emptyset$, 
which contradicts the definition of $\xi_\lambda$.
\end{proof}

\begin{cor} \label{cor:PsubE behaviour} If $\lambda,\mu, \sigma \in \Pi E$ and
$T(\sigma)$ is non-exhaustive, then
\[ 
s_{\sigma\nu_\sigma} s^*_{\sigma\nu_\sigma} s_\lambda s^*_\mu =
\begin{cases} s_{\sigma\nu_\sigma} s^*_{\mu\lambda'\nu_\sigma} 
&\text{if $\sigma =
\lambda\lambda'$} \\ 0 &\text{otherwise.}
\end{cases}
\]
\end{cor}
\begin{proof} The corollary follows from a straightforward calculation using
Lemma~\ref{lem:extensions} and Definition~\ref{dfn:CKfamily}(iii).
\end{proof}

\begin{lemma}\label{lem:QsubE properties} We have
\begin{itemize}
\item[(1)]
$P_{v_0} a \in \lsp\{s_{\lambda\lambda'\nu_{\lambda\lambda'}} s^*_{\mu\lambda'\nu_{\lambda\lambda'}} :
\lambda,\mu \in E, \lambda\lambda' \in\Pi E, 
T(\lambda\lambda')$ non-exhaustive, $s(\nu_{\lambda\lambda'}) = v_0\}$; and
\item[(2)] $\Phi(P_{v_0} a) = P_{v_0} \Phi(a)$.
\end{itemize}
In particular, 
\[\begin{split}
P_{v_0}\Phi(a) \in \lsp\{s_{\lambda\nu_\lambda} s^*_{\mu\nu_\lambda} : {}&\lambda,\mu \in \Pi E,
d(\lambda) = d(\mu), \\
&s(\lambda) = s(\mu), T(\lambda)\text{ non-exhaustive}\}.
\end{split}\]
\end{lemma}
\begin{proof} 
First we use Corollary~\ref{cor:PsubE behaviour} to calculate
\begin{equation} \label{eqn:(1) established} P_{v_0} a  = 
\sum_{\lambda,\mu \in E}  a_{\lambda,\mu}
\bigg(\sum_{\substack{\lambda\lambda' \in \Pi E, 
T(\lambda\lambda')\text{ non-exh.} \\
s(\nu_{\lambda\lambda'}) = v_0}}
s_{\lambda\lambda'\nu_{\lambda\lambda'}} 
s^*_{\mu\lambda'\nu_{\lambda\lambda'}}\bigg)
\end{equation} which proves (1). Furthermore, applying $\Phi$ to \eqref{eqn:(1)
established}, we have
\begin{align*} 
\Phi(P_{v_0} a) &= \sum_{\lambda,\mu \in E} 
a_{\lambda,\mu} \bigg(
                   \sum_{\substack{\lambda\lambda' \in \Pi E, T(\lambda\lambda')\text{ non-exh.} \\
                                   d(\lambda\lambda'\nu_{\lambda\lambda'}) =
d(\mu\lambda'\nu_{\lambda\lambda'}) \\
s(\nu_{\lambda\lambda'}) = v_0}}
                       s_{\lambda\lambda'\nu_{\lambda\lambda'}}
s^*_{\mu\lambda'\nu_{\lambda\lambda'}} \bigg) \\ 
&= \sum_{\substack{\lambda,\mu \in E \\
                    d(\lambda) = d(\mu)}} \bigg(
        a_{\lambda,\mu}
        \sum_{\substack{\lambda\lambda' \in \Pi E, T(\lambda\lambda')\text{ non-exh.} \\
                        s(\nu_{\lambda\lambda'}) = v_0}}
            s_{\lambda\lambda'\nu_{\lambda\lambda'}} 
s^*_{\mu\lambda'\nu_{\lambda\lambda'}}
\bigg) \\ 
&= P_{v_0} \Phi(a).
\end{align*}

The last statement of the lemma follows from (1) and (2) together with 
Remark~\ref{rmk:consequences}(ii).
\end{proof}

We now modify the proof of \cite[Theorem~4.3]{RSY1} to obtain a 
norm-decreasing map
$Q$ which will take $\pi(P_{v_0} a)$ into $\pi(C^*(\L)^\gamma)$.

\begin{lemma} \label{lem:old argument} There exists a norm-decreasing map $Q :
\pi(C^*(\L)) \to \pi(C^*(\L)^\gamma)$ such that
\[
\|Q(\pi(\Phi(P_{v_0}a)))\| = \|\pi(\Phi(P_{v_0}a))\| \text{ and } 
Q(\pi(\Phi(P_{v_0}a)))
= Q(\pi(P_{v_0}a)).
\]
\end{lemma}
\begin{proof} We follow the latter part of the proof of 
\cite[Theorem~4.3]{RSY1} quite
closely. Since $\L$ satisfies \eqref{eqn:aperiodicity}, there exists $x \in
v_0\uL(\infty)$ such that $\lambda \not= \mu$ and $\lambda,\mu \in \L 
v_0$ imply
$\lambda x \not = \mu x$. Hence, for each $\lambda \not= \mu$ in $\L 
v_0$, there exists
$M_{\lambda,\mu} \in \NN^k$ such that $(\lambda x)(0,m) \not= (\mu 
x)(0, m)$ whenever $m
\ge M_{\lambda,\mu}$; assume without loss of generality that 
$M_{\lambda,\mu} \ge
d(\lambda) \vee d(\mu)$. Let
\[\begin{split}
 H:=\{ (\lambda\lambda'\nu_{\lambda\lambda'},\mu\lambda'\nu_{\lambda\lambda'}) :
{}&\lambda,\mu, \lambda\lambda'\in\Pi E, \\
&T(\lambda\lambda')\text{ non-exhaustive}, 
s(\nu_{\lambda\lambda'}) = v_0\},
\end{split}\] 
By Lemma~\ref{lem:QsubE properties}(1),
$P_{v_0} a\in \lsp\{ s_\sigma s^*_\tau :(\sigma,\tau)\in H\}$. Let
\[ T := \{\rho \in \uL({{N_E}}) : \rho = \sigma\text{ or }\rho = 
\tau\text{ for some
}(\sigma,\tau) \in H\}.
\] Define 
\[\textstyle
M:= \bigvee\{M_{\rho,\tau} : \rho \in T, (\sigma,\tau) \in 
H\text{ for some}\sigma,\text{ and } \rho \not= \tau\} + n_x.
\]
The idea is that $M$ is ``far enough
out'' along $x$ to distinguish any pair of paths in $H$. By 
definition of $M$ we have
\begin{equation}\label{eqn:M distingishes} (\tau x)(0, M) \not= (\rho x)(0,M)
\end{equation} when $\tau$ is the second coordinate of an element of 
$H$, $\rho$ belongs
to $T$, and $\tau \not= \rho$. Write $x_M$ for $x(0,M)$.

For $n \le {N_E}$ we set
\[ Q_n := {\sum_{\rho \in T, d(\rho) = n}} \pi(s_{\rho x_M} s^*_{\rho x_M}),
\] and we define $Q : \pi(C^*(\L)) \to \pi(C^*(\L))$ by
\[ Q(b) := \sum_{n \le {N_E}} Q_n b Q_n.
\] As in \cite{RSY1}, $Q$ is norm-decreasing because the $Q_n$ are 
mutually orthogonal
projections. Also as in \cite{RSY1}, $\|Q(\pi(\Phi(P_{v_0}a)))\| =
\|\pi(\Phi(P_{v_0}a))\|$ because $Q$ maps the nonzero matrix units 
in $\pi(P_{v_0} M^s_{\Pi E})$ to nonzero matrix units
in $\pi(\Ff_{N_E + M})$ (see the proof of
\cite[Theorem~4.3]{RSY1} for details).

To establish that $Q(\pi(P_{v_0}a)) = Q(\pi(\Phi(P_{v_0}a)))$, let 
$(\sigma,\tau) \in H$
with $d(\sigma) \not= d(\tau)$. As in the proof of \cite[Theorem~4.3]{RSY1},
$Q(\pi(s_\sigma s^*_\tau))$ is nonzero only if 
there exist $\rho\in T \cap
\L^{d(\sigma)}$ and $\alpha,\beta$ such that
\begin{equation}\label{eqn:paths agree} (\tau x_M \alpha)(0,M) = 
(\rho x_M \beta)(0,M).
\end{equation} 
We claim that $(\tau x_M \alpha)(0,M) = (\tau 
x_M)(0,M)$ for all $\alpha \in s(x_M)\L$: suppose
otherwise for contradiction. Then there exists $i$ such that 
$d(\alpha)_i > 0$ and $d(\tau x_M)_i < M_i$ so $d(x_M)_i <
(M - d(\tau))_i$. But $s((\tau x_M)(0,M)) = s(x_M(0, M-d(\tau)))$, and
since $M \ge d(\tau) + n_x$, we have $M - 
d(\tau) \ge n_x$. It follows that
$\L^{e_i}(x(M - d(\tau))) = \emptyset$ by \eqref{eqn:bdry path}.
The factorisation property now gives
$s(x_M)\L^{e_i} = \emptyset$, contradicting $d(\alpha)_i > 0$. The 
same argument gives
$(\rho x_M \beta)(0,M) = (\rho x_M)(0,M)$ for all $\beta$. 
So \eqref{eqn:paths agree} is equivalent to
$(\tau x_M)(0,M) = (\rho x_M)(0,M)$ which is impossible by 
\eqref{eqn:M distingishes}. Hence $Q(\pi(s_\sigma s^*_\tau)) = 0$
as required.
\end{proof}

\begin{proof}[Proof of Theorem~\ref{thm:CKut}] By \eqref{eqn:norm 
equality}, we have
$\|\Phi(a)\| = \|P_{v_0} \Phi(a)\|$, and Lemma~\ref{lem:QsubE properties} gives
\[\begin{split}
P_{v_0} \Phi(a) \in \lsp\{s_{\lambda\nu_\lambda} 
s^*_{\mu\nu_{\lambda}} : {}&\lambda,\mu\in \Pi E, d(\lambda) = d(\mu) \\
&s(\lambda) = s(\mu), T(\lambda)\text{ non-exhaustive}\}.
\end{split}\] 
Since $\pi$ is injective on the core by 
Theorem~\ref{thm:faithful on core}, we
therefore have
\begin{equation} \label{eqn:QsubE preserves norm}
\|\pi(\Phi(a))\| = \|\Phi(a)\| = \|P_{v_0}\Phi(a)\|= \|\pi(P_{v_0} \Phi(a))\|.
\end{equation} Using \eqref{eqn:QsubE preserves norm},
Lemma~\ref{lem:QsubE properties}(2), and Lemma~\ref{lem:old argument}, 
we therefore have
\begin{align*}
\|\pi(\Phi(a))\| &= \|\pi(P_{v_0} \Phi(a))\| = \|\pi(\Phi(P_{v_0} a))\| \\
&= \|Q(\pi(\Phi(P_{v_0} a)))\| =\|Q(\pi(P_{v_0} a))\| \\ 
&\le \|\pi(P_{v_0}) \pi(a)\| \le \|\pi(a)\|.
\end{align*} The result then follows from 
Proposition~\ref{prp:standard argument}.
\end{proof}

\appendix
\section{The Cuntz-Krieger relations} \label{app:why new rel} The 
objective of the
Cuntz-Krieger relations is to associate to each finitely aligned 
$k$-graph $\L$ a
universal $C^*$-algebra $C^*(\L)$ generated by partial isometries 
$\{s_\lambda :
\lambda \in \L\}$ which has the following properties:
\begin{itemize}
\item[(a)] The partial isometries $s_\lambda$ are all nonzero.
\item[(b)] Connectivity in $\L$ is modelled by multiplication in $C^*(\L)$.
\item[(c)]$C^*(\L)$ is spanned by the elements $\{s_\lambda s^*_\mu : 
\lambda,\mu \in
\L\}$.
\item[(d)] The \emph{core} subalgebra $\clsp\{s_\lambda s^*_\mu : 
\lambda,\mu \in \L,
d(\lambda) = d(\mu)\}$ is AF.
\item[(e)] A representation $\pi$ of $C^*(\L)$ is faithful on the 
core if and only if
$\pi(s_v) \not= 0$ for every vertex $v$.
\end{itemize}

Relations (i) and (ii) of Definition~\ref{dfn:CKfamily} address property (b).
Definition~\ref{dfn:CKfamily}(iii) ensures that property (c) is satisfied.
Definition~\ref{dfn:CKfamily}(iii) has not appeared explicitly in 
previous analyses of
Cuntz-Krieger algebras, but it has always been a consequence of the 
Cuntz-Krieger
relations (see, for example, \cite[Proposition~3.5]{RSY1}). Proposition~6.4 of
\cite{RS1} indicates why we have to impose 
Definition~\ref{dfn:CKfamily}(iii) explicitly
to deal with $k$-graphs that are not row-finite. The analysis of
Section~\ref{sec:thecore} shows that relations (i)--(iii) of
Definition~\ref{dfn:CKfamily} also guarantee property (d).

We must now produce a fourth Cuntz-Krieger relation which guarantees 
that $C^*(\L)$ satisfies (a) and (e); in the following discussion, 
therefore, we assume that Definition~\ref{dfn:CKfamily}(i)--(iii) hold.
We describe examples of $k$-graphs using their \emph{$1$-skeletons\/} as in
\cite[Section~2]{RSY1}. 

The analyses of \cite{FLR} and \cite{RSY1} suggest that a suitable 
relation might be
\begin{equation}\label{eqn:wrong rel 1}
\textstyle t_v = \sum_{\lambda \in v\uL(n)} t_\lambda 
t^*_\lambda\text{ whenever }
v\uL(n)\text{ is finite.}
\end{equation} 
However, this relation fails to guarantee (a), 
even for row-finite $k$-graphs, as 
can be seen from the following example:

\begin{example} \label{ex:non-loc-conv} Consider the row-finite 
2-graph $\L_1$ with
1-skeleton
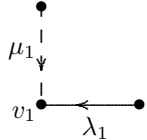
\[
\beginpicture
\setcoordinatesystem units <1.3cm,1.3cm>
\put{$\bullet$} at 1 -.5
\put{$\bullet$} at 2 -.5
\put{$\bullet$} at 1 .5
\setlinear
\setsolid
\plot 1 -.5 2 -.5 /
\arrow <0.15cm> [0.25,0.75] from  1.4 -.5 to 1.36 -.5
\setdashes
\plot 1 -.4 1 .5 /
\arrow <0.15cm> [0.25,0.75] from  1 -.1 to 1 -.14
\put{$\lambda_1$} at 1.55 -.73
\put{$\mu_1$} at .8 .05
\put{$v_1$} at .82 -.59
\endpicture
\] 
where $d(\lambda_1) = (1,0)$ and $d(\mu_1) = (0,1)$. The range 
projections $s_{\lambda_1}
s^*_{\lambda_1}$ and $s_{\mu_1} s^*_{\mu_1}$ are orthogonal by \eqref{eqn:wrong 
rel 1} for $n = (1,1)$, but must both be equal to $s_{v_1}$ by 
\eqref{eqn:wrong rel 1} with $n = (0,1)$ and
$n = (1,0)$. Consequently $s_{v_1} = 0$, so \eqref{eqn:wrong rel 1} 
fails to ensure condition (a) for $C^*(\L_1)$.
\end{example}

For the \emph{row-finite} $k$-graphs of \cite{RSY1} ($v\L^{e_i}$ is always finite), 
we avoided the problem illustrated by this example by
assuming that our $k$-graphs $(\L,d)$ were \emph{locally convex\/}: 
the $k$-graph $(\L,d)$ is locally convex if
for all $v \in \Lambda^0$, $i \not= j$, $\lambda \in v 
\Lambda^{e_i}$ and
$\mu \in v \Lambda^{e_j}$, both $s(\lambda) \Lambda^{e_j}$ and $s(\mu) 
\Lambda^{e_i}$ are nonempty \cite[Definition~3.9]{RSY1}.

For locally convex row-finite $k$-graphs, 
the Cuntz-Krieger relations used in
\cite{RSY1} are equivalent to Definition~\ref{dfn:CKfamily}(i)--(iii) and
\eqref{eqn:wrong rel 2}. It is shown in \cite[Theorem~3.15]{RSY1} 
that these relations
imply (a), and the discussion of \cite[page~109]{RSY1} shows that they 
imply (e). However,
Example~\ref{ex:non-row-fin} demonstrates that for non-row-finite 
$k$-graphs, local
convexity is not enough to ensure that \eqref{eqn:wrong rel 1} implies (e).

\begin{example}\label{ex:non-row-fin} Consider the locally convex 
finitely aligned
2-graph $\L_2$ with 1-skeleton
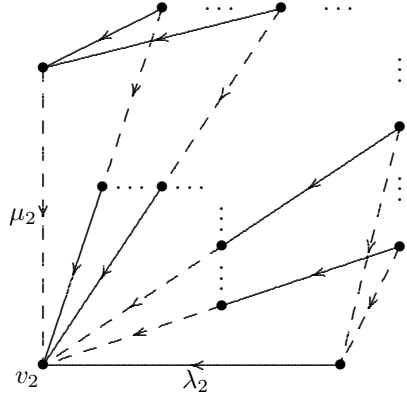
\[
\beginpicture
\setcoordinatesystem units <.75pt,.75pt>
\put{$\bullet$} at 0 0
\put{$v_2$} at -8 -7
\put{$\bullet$} at 150 0
\put{$\bullet$} at 0 150
\put{$\bullet$} at 180 60
\put{$\vdots$} at 180 93
\put{$\bullet$} at 180 120
\put{$\vdots$} at 180 153
\put{$\bullet$} at 60 180
\put{$\dots$} at 90 180
\put{$\bullet$} at 120 180
\put{$\dots$} at 150 180
\put{$\bullet$} at  90 30
\put{$\vdots$} at 90 48
\put{$\bullet$} at 90 60
\put{$\vdots$} at 90 78
\put{$\bullet$} at 30 90
\put{$\dots$} at 45 90
\put{$\bullet$} at 60 90
\put{$\dots$} at 75 90
\plot 150 0 0 0 /
\plot 60 180 0 150 /
\plot 120 180 0 150 /
\plot 30 90 0 0 /
\plot 60 90 0 0 /
\plot 180 60 90 30 /
\plot 180 120 90 60 /
\arrow <0.15cm> [0.25,0.75] from 75.02 0 to 74.98 0
\put{$\lambda_2$} at 77.0 -8.0
\arrow <0.15cm> [0.25,0.75] from 29.0 164.5 to 28.6 164.3
\arrow <0.15cm> [0.25,0.75] from 58.0 164.5 to 57.2 164.3
\arrow <0.15cm> [0.25,0.75] from 15.5 46.5 to 15.3 45.9
\arrow <0.15cm> [0.25,0.75] from 29.8 44.7 to 29.4 44.1
\arrow <0.15cm> [0.25,0.75] from 136.8 45.6 to 136.2 45.4
\arrow <0.15cm> [0.25,0.75] from 134.7 89.8 to 134.1 89.4
\arrow <0.15cm> [0.25,0.75] from 0 75.02 to 0 74.98
\put{$\mu_2$} at -10.0 75.0
\arrow <0.15cm> [0.25,0.75] from 164.5 29.0 to 164.3 28.6
\arrow <0.15cm> [0.25,0.75] from 164.5 58.0 to 164.3 57.2
\arrow <0.15cm> [0.25,0.75] from 46.5 15.5 to 45.9 15.3
\arrow <0.15cm> [0.25,0.75] from 44.7 29.8 to 44.1 29.4
\arrow <0.15cm> [0.25,0.75] from 45.6 136.8 to 45.4 136.2
\arrow <0.15cm> [0.25,0.75] from 89.8 134.7 to 89.4 134.1
\setdashes
\plot 0 150 0 0 /
\plot 180 60 150 0 /
\plot 180 120 150 0 /
\plot 90 30 0 0 /
\plot 90 60 0 0 /
\plot 60 180 30 90 /
\plot 120 180 60 90 /
\endpicture
\] where solid edges have degree $(1,0)$ and dashed edges have degree $(0,1)$.
Relation~\eqref{eqn:wrong rel 1} does not impose any equalities at 
$v_2$ because $v_2 \L^{\le n}_2$ is infinite for all $n \not= 0$. 
The Cuntz-Krieger family $\{S_\lambda :
\lambda \in \L_2\}$ provided by the boundary-path representation 
satisfies $S_{v_2} -
(S_{\lambda_2} S^*_{\lambda_2} + S_{\mu_2} S^*_{\mu_2}) = 0$. 
However, for any nontrivial projection $P$, 
taking $T_{v_2} := S_{v_2} \oplus P$ and 
$T_\sigma = S_\sigma \oplus 0$ for $\sigma \in \L_2 \setminus\{v_2\}$
gives a Cuntz-Krieger $\L_2$-family satisfying 
Definition~\ref{dfn:CKfamily}
(i)--(iii) and \eqref{eqn:wrong rel 1} in which $T_{v_2} - (T_{\lambda_2}
T^*_{\lambda_2} + T_{\mu_2}
T^*_{\mu_2}) \not= 0$. In particular, 
$\{S_\lambda : \lambda \in \L^2\}$ satisfies
Definition~\ref{dfn:CKfamily} (i)--(iii) and \eqref{eqn:wrong rel 1}, but the
representation determined by 
$\{S_\lambda : \lambda \in \L^2\}$ is not faithful on the 
core, even though $S_v \not= 0$ for all $v \in \L_2^0$.
\end{example}

The key property of $\L_2$ which causes the problems with 
relation~\eqref{eqn:wrong rel 1} is
that there exists a finite subset of $v_2 \L_2$ (namely 
$\{\lambda_2,\mu_2\}$) whose range
projections together dominate all the range projections associated to 
paths in $v_2 \L_2 \setminus\{v\}$, 
but no such subset of the form 
$v_2\L^{\le n}_2$. For a
finitely aligned $k$-graph $\L$ and $v \in \L^0$, we can use
Definition~\ref{dfn:CKfamily}(iii) to characterise the finite subsets 
of $v\L$ whose range projections together dominate 
all the range projections associated to nontrivial
paths with range $v$: they are precisely the finite exhaustive sets of
Definition~\ref{dfn:exhaustive}.

Example~\ref{ex:non-row-fin} therefore suggests that Cuntz-Krieger relation 
(iv) should be
\begin{equation}\label{eqn:wrong rel 2}
\textstyle t_v = \sum_{\lambda \in E} t_\lambda 
t^*_\lambda\quad\text{for every $v \in \L^0$ and 
finite exhaustive $E \subset v\L \setminus\{v\}$.}
\end{equation}

\begin{example*}[Example~\ref{ex:non-loc-conv} continued] The only 
finite exhaustive
subset of $v_1 \L_1$ which does not contain $v_1$ is the set 
$\{\lambda_1,\mu_1\}$. In
particular, \eqref{eqn:wrong rel 2} does not insist that either 
$t_{\lambda_1} t_{\lambda_1}^*$
or $t_{\mu_1} t^*_{\mu_1}$ is equal to $t_{v_1}$, and so replacing 
\eqref{eqn:wrong rel 1} with
\eqref{eqn:wrong rel 2} eliminates the pathology associated to the 
non-local-convexity of $\L_1$.
\end{example*}

The only problem with \eqref{eqn:wrong rel 2} is that it is 
predicated on the notion
that the range projections associated to paths in a finite exhaustive 
subset of $v \L \setminus \{v\}$ are mutually orthogonal. The following 
example shows that this is not true.

\begin{example}\label{ex:worst case} Consider the locally convex 
2-graph $\L_3$ with
1-skeleton
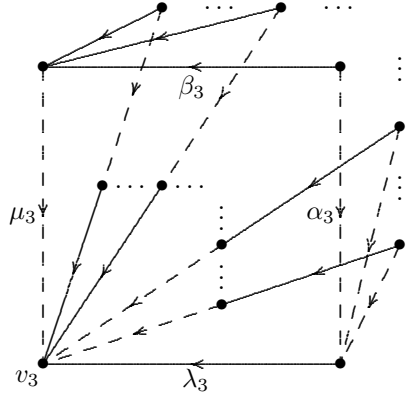
\[
\beginpicture
\setcoordinatesystem units <.75pt,.75pt>
\put{$\bullet$} at 0 0
\put{$v_3$} at -8 -7
\put{$\bullet$} at 150 0
\put{$\bullet$} at 0 150
\put{$\bullet$} at 180 60
\put{$\vdots$} at 180 93
\put{$\bullet$} at 180 120
\put{$\vdots$} at 180 153
\put{$\bullet$} at 60 180
\put{$\dots$} at 90 180
\put{$\bullet$} at 120 180
\put{$\dots$} at 150 180
\put{$\bullet$} at  90 30
\put{$\vdots$} at 90 48
\put{$\bullet$} at 90 60
\put{$\vdots$} at 90 78
\put{$\bullet$} at 30 90
\put{$\dots$} at 45 90
\put{$\bullet$} at 60 90
\put{$\dots$} at 75 90
\put{$\bullet$} at 150 150
\setdashes
\plot 150 150 150 0 /
\setsolid
\plot 150 150 0 150 /
\arrow <0.15cm> [0.25,0.75] from 150 76 to 150 74
\arrow <0.15cm> [0.25,0.75] from 76 150 to 74 150
\put{$\alpha_3$} at 140 75
\put{$\beta_3$} at 75 140
\plot 150 0 0 0 /
\plot 60 180 0 150 /
\plot 120 180 0 150 /
\plot 30 90 0 0 /
\plot 60 90 0 0 /
\plot 180 60 90 30 /
\plot 180 120 90 60 /
\arrow <0.15cm> [0.25,0.75] from 75.02 0 to 74.98 0
\put{$\lambda_3$} at 77.0 -8.0
\arrow <0.15cm> [0.25,0.75] from 29.0 164.5 to 28.6 164.3
\arrow <0.15cm> [0.25,0.75] from 58.0 164.5 to 57.2 164.3
\arrow <0.15cm> [0.25,0.75] from 15.5 46.5 to 15.3 45.9
\arrow <0.15cm> [0.25,0.75] from 29.8 44.7 to 29.4 44.1
\arrow <0.15cm> [0.25,0.75] from 136.8 45.6 to 136.2 45.4
\arrow <0.15cm> [0.25,0.75] from 134.7 89.8 to 134.1 89.4
\arrow <0.15cm> [0.25,0.75] from 0 75.02 to 0 74.98
\put{$\mu_3$} at -10.0 75.0
\arrow <0.15cm> [0.25,0.75] from 164.5 29.0 to 164.3 28.6
\arrow <0.15cm> [0.25,0.75] from 164.5 58.0 to 164.3 57.2
\arrow <0.15cm> [0.25,0.75] from 46.5 15.5 to 45.9 15.3
\arrow <0.15cm> [0.25,0.75] from 44.7 29.8 to 44.1 29.4
\arrow <0.15cm> [0.25,0.75] from 45.6 136.8 to 45.4 136.2
\arrow <0.15cm> [0.25,0.75] from 89.8 134.7 to 89.4 134.1
\setdashes
\plot 0 150 0 0 /
\plot 180 60 150 0 /
\plot 180 120 150 0 /
\plot 90 30 0 0 /
\plot 90 60 0 0 /
\plot 60 180 30 90 /
\plot 120 180 60 90 /
\endpicture
\] where solid edges have degree $(1,0)$ and dashed edges have degree 
$(0,1)$. As in
Example~\ref{ex:non-row-fin}, the fourth Cuntz-Krieger relation must 
insist that the
range projections associated to $\lambda_3$ and $\mu_3$ together fill up 
$t_{v_3}$, or else (e)
will fail because $\{\lambda_3,\mu_3\}$ is finite and exhaustive. 
However, the range
projections $t_{\lambda_3} t^*_{\lambda_3}$ and 
$t_{\mu_3} t^*_{\mu_3}$ are not orthogonal: by
Lemma~\ref{lem:uptoLambda orth range projections}(i), 
$t_{\lambda_3} t^*_{\lambda_3} t_{\mu_3} t^*_{\mu_3} 
= t_{\lambda_3\alpha_3} t^*_{\lambda_3\alpha_3}$. Indeed there is no 
finite exhaustive subset of $v \L$ whose range projections are orthogonal.
\end{example}

The solution to the problem illustrated in Example~\ref{ex:worst 
case} is to use
products rather than sums to express the fourth Cuntz-Krieger relation.

\begin{example*}[Example~\ref{ex:worst case} continued] 
Lemma~\ref{lem:uptoLambda orth
range projections}(i) says that in any family satisfying
Definition~\ref{dfn:CKfamily}(i)--(iii), the projections 
$t_{\lambda_3} t^*_{\lambda_3}$ and
$t_{\mu_3} t^*_{\mu_3}$ commute. 
Consequently, it makes sense to express the 
requirement that
the range projections associated to $\lambda_3$ and $\mu_3$ fill up $t_{v_3}$ with the formula
\begin{equation}\label{eqn:specific rel} (t_{v_3} - t_{\lambda_3} 
t^*_{\lambda_3})(t_{v_3} - t_{\mu_3}
t^*_{\mu_3}) = 0.
\end{equation}
\end{example*}

Relation (iv) of Definition~\ref{dfn:CKfamily}, namely
\begin{equation}\label{eqn:right rel}
\begin{split}
\textstyle \prod_{\lambda \in E} (t_v -  t_\lambda t^*_\lambda) = 0 
\quad\text{for every $v \in \L^0$ and finite exhaustive $E \subset v\L$,}
\end{split}
\end{equation}  is the generalisation of \eqref{eqn:specific rel} to 
arbitrary finite
exhaustive sets in an arbitrary finitely aligned $k$-graph. Note that
\eqref{eqn:right rel} reduces to \eqref{eqn:wrong rel 2} when the 
range projections
associated to paths in $E$ are mutually orthogonal (as in $\L_2$). 
Proposition~\ref{prp:bdry repn
nonzero} together with Theorem~\ref{thm:faithful on core} show that 
\eqref{eqn:right
rel} ensures (a) and (e).

\section{$1$-graphs and locally convex row-finite $k$-graphs} 
\label{app:new rels
generalise}

Recall from \cite{RSY1} that a $k$-graph $(\L,d)$ is \emph{row-finite\/} if
$v \L^{e_i}$ is finite for all $i \in \{1,\dots,k\}$ and $v \in \L^0$. Recall also from
\cite{RSY1} that $(\L,d)$  is \emph{locally convex} if $\lambda \in v \L^{e_i}$ and 
$v\L^{e_j} \not=\emptyset$ for $i \not= j$ implies $s(\lambda) \L^{e_j} \not= 
\emptyset$.

\begin{prop}\label{prop:new rels generalise} For 1-graphs, the 
Cuntz-Krieger families of
Definition~\ref{dfn:CKfamily} coincide with those of \cite{FLR}. For 
locally convex
row-finite $k$-graphs, the  the Cuntz-Krieger families of 
Definition~\ref{dfn:CKfamily}
coincide with those of \cite{RSY1}.
\end{prop}

We prove Proposition~\ref{prop:new rels generalise} with three Lemmas.

\begin{lemma} \label{lem:new families generalise} Let $(\L,d)$ be a 
$k$-graph. If $k >1$, suppose that $\Lambda$ is locally convex and row-finite. Let 
$\{t_\lambda : \lambda
\in \L\}$ be a Cuntz-Krieger $\L$-family. Then $\{t_\lambda : \lambda 
\in \Lambda\}$ is
a Cuntz-Krieger $\L$-family in the sense of \cite{FLR} if $k = 1$, and is a
Cuntz-Krieger $\L$-family in the sense of \cite{RSY1} if $k > 1$.
\end{lemma}

\begin{proof} By Lemma~\ref{lem:uptoLambda orth range 
projections}(iii), we know that
$t_v \ge \sum_{\lambda \in E} t_\lambda t^*_\lambda$ whenever 
$E \subset v\L^{e_i}$ is finite.
By \cite[Propostion~3.11]{RSY1}, it suffices to show that for every $v \in \L^0$ and $1 \le i 
\le k$ such that $0 <
|v \L^{e_i}| < \infty$, we have
\[
 t_v = {\sum_{\lambda \in v \L^{e_i}}} t_\lambda t^*_\lambda.
\] 
By Definition~\ref{dfn:CKfamily}(iv), we need only show that $v \L^{e_i}$ is
exhaustive whenever $0 < |v \L^{e_i}| < \infty$. This is trivial for 
$k = 1$: every path
with range $v$ is either equal to $v$, in which case it is extended 
by every path in $v
\L^{e_1}$, or has an initial segment of length 1, and hence must 
extend an edge in
$\L^{e_1}$. Now suppose $k > 1$ and $\L$ is locally convex and row-finite, 
fix $v$, $i$ with $v \L^{e_i} \not= 
\emptyset$, and let
$\lambda \in v\L$. We must show that there exists $\mu\in v\L^{e_i}$ such that
$\Lmin(\lambda,\mu)\neq\emptyset$ .  If $\lambda = v$, then 
$\Lmin(\lambda,\mu) =
\{(\mu,s(\mu))\}$ for all $\mu \in v \L^{e_i}$. If $d(\lambda) \ge e_i$, then
with $\mu = \lambda(0,e_i) \in v\L^{e_i}$, we have $\Lmin(\lambda,\mu) = 
\{(s(\lambda), \lambda(e_i, d(\lambda)))\} \not= \emptyset$.
Finally, if $\lambda \not= v$ and $d(\lambda)_i = 0$,
then since $v\L^{e_i}$ is nonempty, $|d(\lambda)|$ applications of the local 
convexity condition show that there exists $\alpha \in s(\lambda)\L^{e_i}$.
With $\mu := (\lambda\alpha)(0, e_i)$ 
and $\beta:= (\lambda\alpha)(e_i, d(\lambda\alpha))$ we
have $\mu \in v \L^{e_i}$ and $(\alpha,\beta) \in \Lmin(\lambda,\mu)$.
\end{proof}

\begin{lemma} \label{lem:CKiv FLR} Let $\L$ be a $1$-graph and 
suppose that $\{t_\lambda
: \lambda \in \L\}$ is a Cuntz-Krieger $\L$-family in the sense of 
\cite{FLR}. Then
$\{t_\lambda : \lambda \in \Lambda\}$ satisfies {\rm(iv)} of
Definition~\ref{dfn:CKfamily}.
\end{lemma}

\begin{proof} Let $v\in\L^0$ and let $E$ be a finite exhaustive 
subset of $v\L$. We
proceed by induction on $L(E)  := |\{i \in \NN : E \cap \L^i \not= 
\emptyset\}|$. For a
basis case, suppose that $L(E) = 1$, so $E \subset \L^i$ for some 
$i$. Then $\{\lambda(0,j) : \lambda \in E\} 
= v \L^j$ for $1
\le j \le i$, and then $i$ applications of \cite[Equation~(1.3)]{FLR} give
\[
\prod_{\lambda \in E} (s_v - s_\lambda s^*_\lambda) = s_v - 
\sum_{\lambda \in E} s_\lambda s^*_\lambda = 0.
\]

Now fix $l \ge 1$ and suppose that Definition~\ref{dfn:CKfamily}(iv) holds whenever
$L(E) \le l$, and suppose that $L(E) = l+1$.  Let $I := \max\{i : E \cap \Lambda^i \not= 
\emptyset\}$. Since $L(E) \ge
2$, $\{\lambda \in E : d(\lambda) < I\}$ is nonempty, so let $J := 
\max\{j < I : E \cap
\Lambda^j \not = \emptyset\}$. Fix $\lambda \in E$ with $d(\lambda) = 
I$. Since $E$ is
exhaustive, we have either $\lambda(0,j) \in E$ for some $j \le J$
or $\{\lambda(0,J)\nu : \nu \in
s(\lambda(0,J))\L^{I-J}\} \subset E$. If $\lambda(0,j) \in E$ for some $j \le J$, 
then $t_v - t_\lambda t^*_\lambda \ge t_v - t_{\lambda(0,j)} t^*_{\lambda(0,j)}$, and $E' := E
\setminus\{\lambda\}$ is exhaustive with $\prod_{\mu \in E'} (s_v - 
s_\mu s^*_\mu) = \prod_{\mu \in E} (s_v - s_\mu s^*_\mu)$. On the other hand, if 
$\{\lambda(0,J)\nu : \nu
\in s(\lambda(0,J))\L^{I-J}\} \subset E$, then
\[ E' := \big(E\setminus\{\lambda(0,J)\nu : \nu \in 
s(\lambda(0,J))\L^{I-J}\}\big)
\cup\{\lambda(0,J)\}
\] is also exhaustive, and $\prod_{\mu \in E'} (s_v - s_\mu s^*_\mu) 
= \prod_{\mu \in E}
(s_v - s_\mu s^*_\mu)$. Repeating this process for each $\lambda \in 
E \cap \L^I$, we
obtain a finite exhaustive $E'' \in v\L$ which satisfies
\begin{itemize}
\item[(1)] $\{i \in \NN : E'' \cap \Lambda^i \not= \emptyset\} = \{i 
\in \NN : E \cap
\Lambda^i \not= \emptyset\} \setminus\{ I \}$, so $L(E'') = L(E) - 1 = l$; and
\item[(2)] $\prod_{\mu \in E''} (s_v - s_\mu s^*_\mu) = \prod_{\mu 
\in E} (s_v - s_\mu
s^*_\mu)$.
\end{itemize} The result now follows from the inductive hypothesis 
applied to $E''$.
\end{proof}

\begin{lemma} \label{lem:CKiv RSY1} Let $(\L,d)$ be a locally convex row-finite
$k$-graph and let $\{t_\lambda : \lambda \in \Lambda\}$ be a 
Cuntz-Krieger $\L$-family
in the sense of \cite[Definition~3.3]{RSY1}.  Then $\{t_\lambda : 
\lambda \in \Lambda\}$
satisfies {\rm(iv)} of Definition~\ref{dfn:CKfamily}.
\end{lemma}

\begin{proof} Let $v\in\L^0$, let $E$ be a finite exhaustive subset 
of $v\L$, and let $N:= \bigvee_{\lambda \in E} d(\lambda)$. 
Now let $E' := \{\lambda\nu : \lambda \in E,
\nu\in s(\lambda)\uL({N-d(\lambda)})\}$. By \cite[Lemma~3.6]{RSY1}, 
and since $E$ is exhaustive, we have $E' = v\uL(N)$.
Hence relation (4) of \cite[Definition~3.3]{RSY1} ensures that 
$s_v = \sum_{\mu \in E'} s_\mu s^*_\mu$,
so
\begin{flalign*}
&&\prod_{\lambda \in E} (s_v - s_\lambda s^*_\lambda) 
&\le \prod_{\mu \in E'} (s_v - s_\mu s^*_\mu)
= s_v - \sum_{\mu \in v\uL(N)} s_\mu s^*_\mu
=0. &\qed
\end{flalign*}
\renewcommand\qed{}\end{proof}

\begin{proof}[Proof of Proposition~\ref{prop:new rels generalise}]
 Lemma~\ref{lem:new families generalise} shows that the 
Cuntz-Krieger
families of Definition~\ref{dfn:CKfamily} give Cuntz-Krieger families
as defined in \cite{FLR} and
\cite{RSY1}. Relations (i) and (ii) of Definition~\ref{dfn:CKfamily} 
are obviously
satisfied by the Cuntz-Krieger families of both \cite{FLR} and 
\cite{RSY1}.  In a
1-graph, $\Lmin(\lambda,\mu)$ equals $\{(\lambda',s(\mu))\}$ if $\mu 
= \lambda\lambda'$,
$\{(s(\lambda),\mu')\}$ if $\lambda = \mu\mu'$, and $\emptyset$ 
otherwise. It follows that
relation~(iii) of Definition~\ref{dfn:CKfamily} is satisfied by the 
Cuntz-Krieger families of \cite{FLR}. Proposition~3.5 of \cite{RSY1} shows that for 
locally convex row-finite $k$-graphs, Relation~(iii) of 
Definition~\ref{dfn:CKfamily} is satisfied by
the Cuntz-Krieger families of \cite{RSY1}. The result now follows from
Lemmas~\ref{lem:CKiv FLR} and \ref{lem:CKiv RSY1}.
\end{proof}

\section{Checking the relations in terms of 
generators}\label{app:gens do} 

\begin{theorem} \label{thm:gen families}
Let $(\L,d)$ be a finitely aligned $k$-graph. Let 
\[\textstyle
\big\{t_\lambda : \lambda \in \big(\bigcup^k_{i=1} \L^{e_i}\big) 
\cup \L^0\big\}
\]
be a family of partial isometries in a $C^*$-algebra. Then there is at most one Cuntz-Krieger $\L$-family
$\{t'_\lambda : \lambda \in \L\}$ such that $t'_\lambda = t_\lambda$ for all 
$\lambda \in \big(\bigcup^k_{i=1} \L^{e_i}\big) \cup \L^0\big)$. Furthermore, such a Cuntz-Krieger
$\L$-family exists if and only if
\begin{itemize}
\item[(i)] $\{t_v : v \in \L^0\}$ is a collection of mutually orthogonal projections.
\item[(ii)] $t_\lambda t_\alpha = t_\mu t_\beta$ when 
$\lambda, \mu, \alpha, \beta \in \big(\bigcup^k_{i=1} \L^{e_i}\big) \cup \L^0$ 
satisfy $\lambda\alpha = \mu\beta$. 
\item[(iii)] $t^*_\lambda t_\mu = \sum_{(\alpha,\beta) \in \Lmin(\lambda,\mu)} t_\alpha t^*_\beta$ for all $\lambda,\mu \in \bigcup^k_{i=1} \L^{e_i}$.
\item[(iv)] for every $v \in \L^0$ and every finite exhaustive $E \subset 
\bigcup^k_{i=1} v\L^{e_i}$,
\[
\prod_{\lambda \in E} (t_v - t_\lambda t^*_\lambda) = 0.
\]
\end{itemize}
\end{theorem}

Before proving Theorem~\ref{thm:gen families}, we establish a number of preliminary results.

\begin{lemma}\label{lem:CKiii gens}
 Let $(\L,d)$ be a finitely aligned $k$-graph. Suppose 
that $\{t_\lambda :
\lambda \in \L\}$ is a collection of partial isometries satisfying
Definition~\ref{dfn:CKfamily}(i) and (ii). Then $\{t_\lambda : 
\lambda \in \L\}$
satisfies  Definition~\ref{dfn:CKfamily}(iii) if and only if
\begin{equation}\label{eqn:gens(iii)}
\textstyle t^*_\lambda t_\mu = \sum_{(\alpha,\beta) \in 
\Lmin(\lambda,\mu)} t_\alpha
t^*_\beta \text{ for all } \lambda, \mu \in \bigcup^k_{i=1} \L^{e_i}.
\end{equation}
\end{lemma}
\begin{proof} Since \eqref{eqn:gens(iii)} is a special case of
Definition~\ref{dfn:CKfamily}(iii), we need only show the ``if'' 
direction. This in
turn will follow from \cite[Lemma~9.2]{RS1} if we can show that
Definition~\ref{dfn:CKfamily}(i) and (ii) together with \eqref{eqn:gens(iii)} 
imply
relations (3) and (4) of \cite[Definition~7.1]{RS1}, namely that
\begin{gather} t^*_\lambda t_\lambda = t_{s(\lambda)}\text{ for all 
$\lambda \in \L$;
and} \label{eqn:oldrel3} \\
\textstyle t_v \ge \sum_{\lambda \in F} t_\lambda t^*_\lambda\text{ 
whenever $F \subset
\L^n v$ is finite.}
\label{eqn:oldrel4}
\end{gather} An inductive argument on the length of $\lambda$
establishes \eqref{eqn:oldrel3}. With this in hand, 
\eqref{eqn:oldrel4} then follows
from \eqref{eqn:gens(iii)} together with Definition~\ref{dfn:CKfamily}(ii) as in
Lemma~\ref{lem:uptoLambda orth range projections}(iii).
\end{proof}

\begin{prop} \label{prp:gens do} Let $(\L,d)$ be a finitely 
aligned $k$-graph. A
family  $\{t_\lambda : \lambda \in \L\}$ of partial isometries satisfying
Definition~\ref{dfn:CKfamily}(i)--(iii)  is a Cuntz-Krieger 
$\L$-family if and only if
for every $v \in \L^0$ and every finite exhaustive subset $E \subset 
\bigcup^k_{i=1}
v\L^{e_i}$,
\begin{equation}\label{eqn:gen relation}
\prod_{\lambda \in E} (t_v - t_\lambda t^*_\lambda) = 0.
\end{equation}
\end{prop}

\begin{notation} In this section, we make use of the 
following notation:
\begin{itemize}
\item Given a set $E \subset \L$, define
$I(E) := \bigcup^k_{i=1} \{\lambda(0,e_i) : \lambda \in E, d(\lambda)_i > 
0\}$.
\item Given $E \subset \L$ and $\mu \in \L$, let $\Ext(\mu;E) :=\bigcup_{\lambda \in E} \{\alpha : (\alpha,\beta) \in
\Lmin(\mu,\lambda)\}$.
\item Given $E \subset \L$, let $L(E)
:= \sum^k_{i=1} \max_{\lambda \in E} d(\lambda)_i$.
\end{itemize}
\end{notation}

\begin{lemma} \label{lem:exhaustive sets} Let $(\L,d)$ be a finitely aligned $k$-graph
and let $v \in \L^0$. Suppose $E 
\subset v \L$ is finite and exhausitve, and let $\mu \in v \L $. Then
$\Ext(\mu;E)$ is a finite exhaustive subset of $s(\mu) \L$.
\end{lemma}
\begin{proof} Since $E$ is finite and $\L$ is finitely aligned we 
know that $\Ext(\mu;E)$
is finite, so we need only check that $\Ext(\mu;E)$ is exhaustive. 
Let $\sigma \in
s(\mu)\L$. Since $E$ is exhaustive, there exists $\lambda \in E$ with
$\Lmin(\lambda,\mu\sigma)
\not= \emptyset$, say $(\alpha,\beta) \in \Lmin(\lambda,\mu\sigma)$. So
$\lambda\alpha = \mu\sigma\beta$, and hence
\[
\big(\alpha(0,(d(\lambda) \vee d(\mu)) - d(\lambda)), (\sigma\beta) 
(0,(d(\lambda) \vee
d(\mu)) - d(\mu))\big) \in \Lmin(\lambda,\mu).
\] Hence $\tau := (\sigma\beta)(0,(d(\lambda) \vee d(\mu)) - d(\mu))$ 
belongs to
$\Ext(\mu;E)$, and then
\begin{flalign*} &&&\big((\sigma\beta)(d(\sigma),d(\sigma) \vee d(\tau)),
(\sigma\beta)(d(\tau), d(\sigma) \vee d(\tau))\big) \in \Lmin(\sigma,\tau)&\qed
\end{flalign*}
\renewcommand\qed{}\end{proof}

\begin{lemma} \label{lem:e_i exhaustive} Let $(\L,d)$ be a finitely aligned $k$-graph,
let $v \in \L^0$, and 
suppose that $E \subset
v\L \setminus\{v\}$ is finite and exhaustive. Then
$I(E)$ is also finite and exhaustive.
\end{lemma}
\begin{proof} We have $I(E)$ is finite because $E$ is finite, so we 
just need to show
that $I(E)$ is exhaustive. Let $\mu \in v\L$. Since $E$ is 
exhaustive, there exists
$\lambda \in E$ such that $\Lmin(\lambda,\mu) \not= \emptyset$, say 
$(\alpha,\beta)
\in \Lmin(\lambda,\mu)$. Since $\lambda \in E$, we have $d(\lambda) 
\not= 0$, so fix $i$
such that $d(\lambda)_i \not= 0$; then $\lambda(0,e_i) \in I(E)$. Let $\rho :=
(\lambda\alpha)(0, d(\mu) \vee e_i)$, let $\eta := \rho(e_i, 
d(\rho))$, and let $\zeta
:= \rho(d(\mu),d(\rho))$. Then $\lambda(0,e_i) \eta = \rho = 
\mu\zeta$, so $(\eta,\zeta)
\in \Lmin(\lambda(0,e_i), \mu)$. Since $\mu \in v\L$ was arbitrary, 
it follows that
$I(E)$ is exhaustive.
\end{proof}

\begin{lemma} \label{lem:gap projections} Let $(\L,d)$ be a finitely aligned $k$-graph,
and let $\{t_\lambda: \lambda \in \L\}$ be a family of partial isometries satisfying
Definition~\ref{dfn:CKfamily}(i)--(iii). Let $v \in \L^0$, let
$\lambda \in v \L$ and suppose that $E \subset s(\lambda) \L$ is 
finite and satisfies
$\prod_{\nu \in E} (t_{s(\lambda)} - t_\nu t^*_\nu) = 0$. Then
\[ t_v - t_\lambda t^*_\lambda = \prod_{\nu \in E} (t_v - t_{\lambda\nu}
t^*_{\lambda\nu}).
\]
\end{lemma}
\begin{proof} Since $t_{\lambda\mu} t^*_{\lambda\mu} \le t_\lambda 
t^*_\lambda$ for all
$\mu \in s(\lambda)\L$, we have
\[ (t_v - t_\lambda t^*_\lambda)(t_v - t_{\lambda\nu} t^*_{\lambda\nu}) = t_v -
t_\lambda t^*_\lambda
\] for all $\nu \in E$. It follows that
\begin{equation}\label{eqn:LHS} (t_v - t_\lambda t^*_\lambda) 
\prod_{\nu \in E} (t_v -
t_{\lambda\nu} t^*_{\lambda\nu}) = t_v - t_\lambda t^*_\lambda.
\end{equation} On the other hand,
\begin{align*} (t_v - t_\lambda t^*_\lambda) &\Big(\prod_{\nu \in E} (t_v -
t_{\lambda\nu} t^*_{\lambda\nu}) \Big) \\ &= t_v \Big(\prod_{\nu \in E} (t_v -
t_{\lambda\nu} t^*_{\lambda\nu}) \Big) - t_\lambda t^*_\lambda 
\Big(\prod_{\nu \in E}
(t_v - t_{\lambda\nu} t^*_{\lambda\nu})\Big) \\ &= \Big(\prod_{\nu 
\in E} (t_v -
t_{\lambda\nu} t^*_{\lambda\nu})\Big) - \Big(\prod_{\nu \in E} 
(t_\lambda t^*_\lambda -
t_{\lambda\nu} t^*_{\lambda\nu})\Big) \\ &= \Big(\prod_{\nu \in E} 
(t_v - t_{\lambda\nu}
t^*_{\lambda\nu})\Big) - t_\lambda \Big(\prod_{\nu \in E} 
(t_{s(\lambda)} - t_\nu
t^*_\nu)\Big) t^*_\lambda \\ &= \prod_{\nu \in E} (t_v - 
t_{\lambda\nu} t^*_{\lambda\nu})
\end{align*} because $\prod_{\nu \in E} (t_{s(\lambda)} - t_\nu 
t^*_\nu) = 0$ by
hypothesis.
\end{proof}

\begin{lemma} \label{lem:reduces length} Let $(\L,d)$ be a finitely aligned
$k$-graph. Let $v \in \L^0$ and suppose $E \subset v \L$ is finite.
Suppose $\lambda \in I(E)$. Then $L(\Ext(\lambda;E)) < L(E)$.
\end{lemma}
\begin{proof} 
Since $\lambda \in I(E)$, we have $d(\lambda) = e_i$ and $\lambda\lambda' \in E$ for some
$i, \lambda'$. For $j \in \{1, \dots, k\}$, we have
\begin{equation}\label{eqn:maximums}
\max_{\nu \in \Ext(\lambda;E)} d(\nu)_j =
\max_{\mu \in E, \Lmin(\lambda,\mu) \not= \emptyset} ((d(\lambda) \vee 
d(\mu)) - e_i)_j.
\end{equation} If $i \not= j$, then \eqref{eqn:maximums} becomes
\[
\max_{\nu \in \Ext(\lambda;E)} d(\nu)_j =
\max_{\mu \in E, \Lmin(\lambda,\mu) \not= \emptyset} d(\mu)_j
\le \max_{\mu \in E} d(\mu)_j.
\]

On the other hand, if $i = j$, then we use \eqref{eqn:maximums} to calculate
\begin{align*}
\max_{\nu \in \Ext(\lambda;E)} d(\nu)_j 
&= \max_{\mu \in E, \Lmin(\lambda,\mu) \not= 
\emptyset} ((d(\lambda)
\vee d(\mu)) - e_i)_i \\ &\le \max_{\mu \in E} ((d(\lambda) \vee d(\mu)) 
- e_i)_i \\ &=
\big(\max_{\mu \in E} d(\mu)_i\big) - 1 \\ &\quad\text{since $\lambda\lambda' 
\in E$ so there
exist $\mu \in E$ with
$d(\mu)_i \ge 1$}
\end{align*} We therefore have
\begin{flalign*} &&L(\Ext(\lambda;E)) &= \sum_{j=1}^k \max_{\nu \in 
\Ext(\lambda;E)}
d(\nu)_j &\\ &&&\le \Big( \sum_{j \in \{1,\dots,k\}\setminus\{i\}} 
\max_{\mu \in E}
d(\mu)_j \Big) + \big(\max_{\mu \in E} d(\mu)_i \big) - 1 &\\ &&&< 
\sum^k_{j=1} \max_{\mu \in E}
d(\mu)_j &\\ &&&= L(E) &\qed
\end{flalign*}
\renewcommand\qed{}\end{proof}

\begin{proof}[Proof of Proposition~\ref{prp:gens do}] We must show that 
for every $v \in
\L^0$ and every finite exhaustive $F
\subset v\L$, we have
\begin{equation}\label{eqn:needed}
\prod_{\mu \in F} (t_v - t_\mu t^*_\mu)= 0.
\end{equation} We proceed by induction on $L(F)$. If $L(F) = 1$, then 
$F \subset
\bigcup^k_{i=1} v\L^{e_i}$, and \eqref{eqn:needed} is an instance of
\eqref{eqn:gen relation}.

Now suppose that \eqref{eqn:needed} holds whenever $L(F) \le n$, and fix
$v \in \L^0$ and $F \subset v\L$ finite exhaustive with $L(F) = n+1$. If
$v \in F$, there is nothing to prove, so assume without loss of 
generality that $v \not
\in F$. Then $I(F)$ is finite and exhaustive by Lemma~\ref{lem:e_i 
exhaustive}. Fix
$\lambda
\in I(F)$. By Lemma~\ref{lem:exhaustive sets}, we know that 
$\Ext(\lambda;F)$ is finite
and exhaustive. By Lemma~\ref{lem:reduces length}, we know that
$L(\Ext(\lambda;F)) \le n$, so the inductive hypothesis ensures that
$\prod_{\nu \in \Ext(\lambda;F)} (t_{s(\lambda)} - t_\nu t^*_\nu) = 
0$. It then follows
from Lemma~\ref{lem:gap projections} that
\begin{equation}\label{eqn:induct}
\prod_{\nu \in \Ext(\lambda;F)} (t_v - t_{\lambda\nu} t^*_{\lambda\nu}) = t_v -
t_\lambda t^*_\lambda.
\end{equation}

For each $\nu \in \Ext(\lambda;F)$, there exists $\mu \in F$ with $\lambda\nu =
\mu\mu'$, so $t_{\lambda\nu} t^*_{\lambda\nu} \le t_\mu t^*_\mu$, and hence
\begin{equation}\label{eqn:dominates}
\prod_{\nu \in \Ext(\lambda;F)} (t_v - t_{\lambda\nu} 
t^*_{\lambda\nu}) \ge \prod_{\mu
\in F} (t_v - t_\mu t^*_\mu).
\end{equation}

We can therefore calculate
\begin{flalign*} &&\prod_{\mu \in F} (t_v - t_\mu t^*_\mu) &\le 
\prod_{\lambda \in
I(F)}\Big(\prod_{\nu \in \Ext(\lambda;F)} (t_v - t_{\lambda\nu} 
t^*_{\lambda\nu})\Big)
\quad\text{by \eqref{eqn:dominates}} &\\ &&& = \prod_{\lambda \in 
I(F)}(t_v - t_\lambda
t^*_\lambda) \quad\text{by
\eqref{eqn:induct}.} &\\ &&&= 0 \quad\text{by \eqref{eqn:gen relation}}&\qed
\end{flalign*}
\renewcommand\qed{}\end{proof}

\begin{proof}[Proof of Theorem~\ref{thm:gen families}]
The factorisation property and Definition~\ref{dfn:CKfamily}(ii) show that any Cuntz-Krieger $\L$-family
$\{t'_\lambda : \lambda \in \L\}$ satisfying $t'_\lambda = t_\lambda$ for all 
$\lambda \in \big(\bigcup^k_{i=1} \L^{e_i}\big) \cup \L^0$ must satisfy
\begin{equation}\label{eq:tsubmu def} 
t'_\lambda = t_{\lambda_1} t_{\lambda_2} \cdots t_{\lambda_{|d(\lambda)|}}
\end{equation}
for each $\lambda \in \L$ and each factorisation $\lambda = \lambda_1 \cdots \lambda_{|d(\lambda)|}$
where the $\lambda_i$ belong to $\big(\bigcup^k_{i=1} \L^{e_i}\big) \cup \L^0$. This proves that
there is at most one such Cuntz-Krieger $\L$-family.

Suppose that such a Cuntz-Krieger $\L$-family $\{t'_\lambda : \lambda \in \L\}$ exists. Then conditions 
(i)--(iv) of Theorem~\ref{thm:gen families} are immediate consequences of the 
Cuntz-Krieger relations.

Now suppose that 
$\{t_\lambda : \lambda \in \big(\bigcup^k_{i=1} \L^{e_i}\big) \cup \L^0\}$
satisfy (i)--(iv) of Theorem~\ref{thm:gen families}. An inductive argument using
condition~(ii) of Theorem~\ref{thm:gen families} shows that \eqref{eq:tsubmu def}
gives a well-defined family of partial isometries $\{t'_\lambda : \lambda \in \L\}$.

We have that $\{t'_\lambda : \lambda \in \L\}$ satisfies 
Definition~\ref{dfn:CKfamily}(i) because this is precisely 
condition~(i) of Theorem~\ref{thm:gen families}. Equation~\eqref{eq:tsubmu def} 
and the factorisation property for $\L$ ensure that 
$\{t'_\lambda : \lambda \in \L\}$ satisfies Definition~\ref{dfn:CKfamily}(ii).
Condition~(iii) of Theorem~\ref{thm:gen families} and 
Lemma~\ref{lem:CKiii gens} then imply that $\{t'_\lambda : \lambda \in \L\}$ 
satisfies Definition~\ref{dfn:CKfamily}(iii). We can now use 
Proposition~\ref{prp:gens do} and condition~(iv) of 
Theorem~\ref{thm:gen families} to show that $\{t'_\lambda : \lambda \in \L\}$ 
satisfies Definition~\ref{dfn:CKfamily}(iv).
\end{proof}

\end{document}